\documentclass{article}

\usepackage[utf8]{inputenc}
\usepackage[T1]{fontenc}
\usepackage{amsmath,amssymb,amsthm,bbm,float,graphicx,geometry,lmodern,mathtools,parskip,setspace,subcaption}
\usepackage[colorlinks, linkcolor = blue!80!black, citecolor = blue!80!black,breaklinks, pdfauthor={Benedikt Jahnel, Lukas Luechtrath, Anh Duc Vu}]{hyperref}
\usepackage{mathrsfs}
\usepackage{enumerate}
\usepackage{nicefrac}
\usepackage{todonotes}
\usepackage{comment}
\usepackage{titling}
\usepackage{fullpage}
\usepackage{comment}
\usepackage{orcidlink}

\usepackage{tikz}
\usetikzlibrary{backgrounds}
\usetikzlibrary{patterns}
\usetikzlibrary{positioning, shapes.geometric}
\usetikzlibrary{external, arrows, decorations.pathmorphing} 

\usepackage{caption}
\usepackage{graphicx}
\graphicspath{{Images/}}
\allowdisplaybreaks
\usepackage[raggedright]{titlesec} 
\usepackage{booktabs} 

\newcommand{\N}{\mathbb{N}}

\DeclareMathOperator{\interior}{int}

\renewcommand*{\P}{\mathbb{P}}
\newcommand*{\E}{\mathbb{E}}
\newcommand*{\R}{\mathbb{R}}

\newcommand*{\Z}{\mathbb{Z}}

\renewcommand{\N}{\mathbb{N}}
\newcommand*{\X}{\mathbb{X}}

\newcommand*{\Bcal}{\mathcal{B}}

\renewcommand*{\d}{\mathrm{d}}
\newcommand*{\e}{\mathrm{e}}
\renewcommand*{\Xi}{\varXi}
\renewcommand*{\epsilon}{\varepsilon}
\renewcommand*{\theta}{\vartheta}
\renewcommand*{\Theta}{\varTheta}

\renewcommand*{\Delta}{\varDelta}

\newcommand{\g}{\gamma}
\newcommand{\eps}{\varepsilon}

\renewcommand{\interior}[1]{%
  {\kern0pt#1}^{\mathrm{o}}%
}

\newtheorem{theorem}{theorem}[section]
\newtheorem{cor}[theorem]{Corollary}
\newtheorem{lemma}[theorem]{Lemma}
\newtheorem{proposition}[theorem]{Proposition}
\newtheorem{thm}[theorem]{Theorem}
\theoremstyle{definition}

\theoremstyle{remark}
\newtheorem{remark}[theorem]{Remark}
\newtheorem{example}[theorem]{Example}

\usepackage{nameref}

\newcommand{\BBB}{}

\makeatletter
\let\orgdescriptionlabel\descriptionlabel
\renewcommand*{\descriptionlabel}[1]{%
  \let\orglabel\label
  \let\label\@gobble
  \phantomsection
  \edef\@currentlabel{#1}%
  \let\label\orglabel
  \orgdescriptionlabel{#1}%
}

\pretitle{\centering\LARGE\scshape}
 \posttitle{\vskip 0.75cm}

 \predate{\vskip 0.75 cm \centering\large}
 \postdate{\par}

\usepackage[style = numeric, sorting=nyt, url = false, abbreviate=false, maxbibnames=9, sortcites=true, doi = true, backend = biber, giveninits = true, isbn=false]{biblatex}
\renewbibmacro{in:}{\ifentrytype{article}{}{\printtext{\bibstring{in}\intitlepunct}}}
\bibliography{references_2}

\title{First contact percolation}
\thanksmarkseries{arabic}

\author{
Benedikt Jahnel \orcidlink{0000-0002-4212-0065}
\thanksgap{0.4ex}\thanks{Technische Universit\"at Braunschweig, Universit\"atsplatz 2, 38106 Braunschweig, Germany}
\thanksgap{0.4ex}\thanks{Weierstrass Institute for Applied Analysis and Stochastics, Mohrenstraße 39, 10117 Berlin, Germany \\ \phantom{Cor}\(^*\)Corresponding Author: benedikt.jahnel@tu-braunschweig.de} \thanksgap{0.4ex}\(^*\) \\ 
\and
Lukas L\"{u}chtrath \orcidlink{0000-0003-4969-806X} \thanksmark{2}
\and
Anh Duc Vu \orcidlink{0009-0005-6913-4992} \thanksmark{2}
}

\date{August 14, 2025}

\usepackage{cleveref}

\begin{document}

\maketitle

\begin{spacing}{0.9}
\begin{abstract} 
\noindent We study a version of first passage percolation on $\Z^d$ where the random passage times on the edges are replaced by contact times represented by random closed sets on $\R$. Similarly to the contact process without recovery, an infection can spread into the system along increasing sequences of contact times. In case of stationary contact times, we can identify associated first passage percolation models, which in turn establish shape theorems also for first contact percolation. In case of periodic contact times that reflect some reoccurring daily pattern, we also present shape theorems with limiting shapes that are universal with respect to the within-one-day contact distribution. In this case, we also prove a Poisson approximation for increasing numbers of within-one-day contacts. Finally, we present a comparison of the limiting speeds of three models -- all calibrated to have one expected contact per day -- that suggests that less randomness is beneficial for the speed of the infection. The proofs rest on coupling and subergodicity arguments. 

 \medskip
\noindent\footnotesize{{\textbf{AMS-MSC 2020}: Primary: 60K05; Secondary: 60K35, 82B43}

}
\end{abstract}
\end{spacing}

\section{Introduction} \label{sec:intro}
In {\em first passage percolation} (FPP), edges of a graph are assigned random {\em transition times} and one studies the minimal total time needed to connect two vertices, which can be seen as a {\em random metric}~\cite{carmona1986aspects,auffinger2017fpp50years,kesten1987percolation,kesten2003first,smythe2006first}. In particular, one might be interested in the set of vertices that are reached up to some fixed time $t$ from a given vertex (say the origin $o$). This represents, for example, the set of wet sites at time $t$ if $o$ is interpreted as a water source and edges are pipes of random length, justifying the term {\em percolation}. Another common interpretation comes from the modeling of epidemiological processes, where an initially infected network component spreads a virus in a population with varying transition times between neighboring vertices. 

However, one might argue that such a situation is represented more accurately by assigning {\em random contact times} between any pair of vertices connected by an edge. In this case, the time for an infection to reach another vertex is determined by the total time needed along an optimal path of edges with an increasing sequence of contact times; we call this model {\em first contact percolation} (FCP) and dedicate these notes to the study of FCP with respect to the asymptotic speed of infection spread. 

\begin{figure}[t]
\includegraphics[width=1\columnwidth]{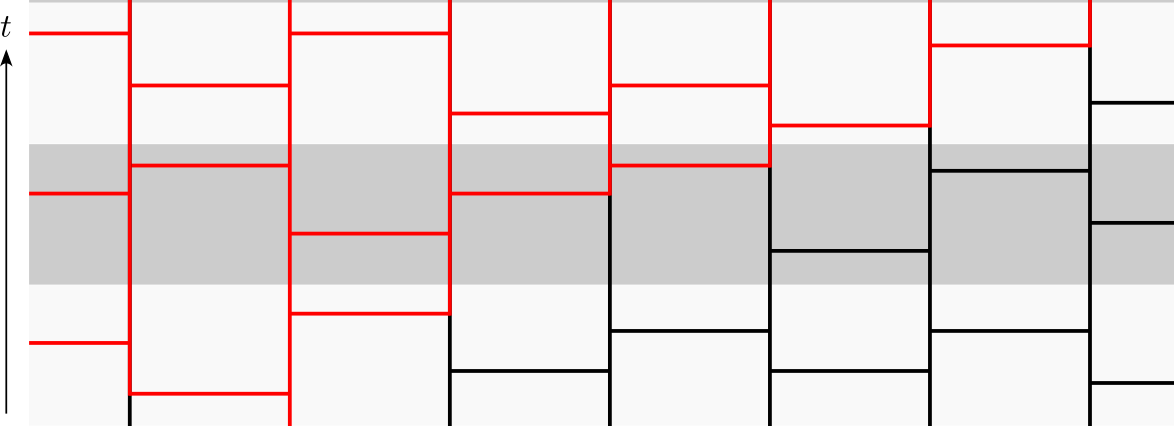}
\caption{Graphical representation of an infection process. Vertical lines indicate vertices and their state. Horizontal lines indicate contact times between neighboring vertices. The infection (red) then spreads across the network along first contact times. In this realization, neighboring vertices have exactly one contact per integer interval (colored background).}\label{fig:graphical-representation}
\end{figure}

It is important to note that the connection times between pairs of vertices $x,y$ no longer form a metric, since the time needed to spread from $x$ to $y$ may differ from the time needed to spread from $y$ to $x$. This asymmetry arises as the traversal time of an edge depends on the time at which one arrives at its endpoint, which is unlike in FPP. The situation is perhaps best illustrated by a graphical representation, see Figure~\ref{fig:graphical-representation}, where each edge is equipped with a point process of contact times and one follows increasing space-time paths. 

\begin{figure}[t]
\includegraphics[width=1\columnwidth]{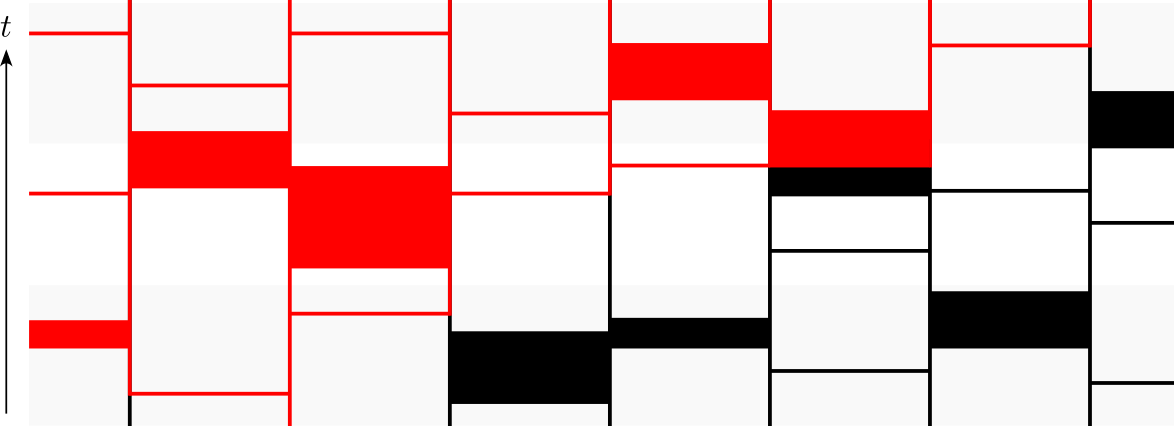}
\caption{Graphical representation of an infection process where contact times can be arbitrary random closed sets, e.g., Boolean models with random radii based on simple point processes.}\label{fig:graphical-representation-boolean}
\end{figure}

An important special case arises when the point processes are i.i.d.\ (with respect to the edges) Poisson point processes. In this case, FCP becomes a Markov process in time, which simply corresponds to the classical {\em contact process} without recovery~\cite{durrett1991contact,durrett1982contact,durrett1980growth}, i.e., its {\em pure-growth} version. In that sense, FCP can be seen as the pure-growth version of contact processes based on general point processes. In fact, one does not have to restrict oneself to point processes on the edges, instead more general {\em random closed sets} of times can be considered, which represent the time and duration of a meeting, see Figure~\ref{fig:graphical-representation-boolean} for an illustration. In this context, let us mention results on the contact process beyond exponential waiting times~\cite{fontes2021contact,fontes2019contact,fontes2020contact,fontes2023renewal,linker2020contact,hilario2022results,seiler2023contact,cardona2024contact}, which primarily focus on either {\em renewal processes} or {\em time-dependent processes} as alternatives to the classical Poisson point processes.  

The main focus of this study of FCP is to derive and compare the asymptotic shape of the infected vertices as time tends to infinity. More precisely, we consider FCP on the hypercubic lattice $\Z^d$ with $d\ge 1$, where each edge $e$ is independently equipped with a random closed set of meeting times $\X_e\subset \R$ following the same distribution. Note that each edge must only be traversed once in order to accurately describe the distribution of the vertices reached up to a given time. This allows us to draw a connection to FPP and ultimately to recover an associated {\em shape theorem} when $\X_e$ is stationary. Such a direct connection is not possible for non-stationary $\X_e$. However, we establish a shape theorem also for specific models that obey stationarity with respect to shifts in $\Z$. Paradigmatically, we consider a situation in which the time axis is split into {\em days} $\{[x,x+1)\colon x\in \Z\}$ and, \(n\in\N\) i.i.d.\ meeting times are assigned per day and edge. 
We observe that the limiting shape of infected vertices in this model is independent of the underlying diffusive distribution of meeting times since only order statistics are relevant, which allows for the restriction to uniformly distributed meeting times. 
Furthermore, the properly rescaled limiting shape obeys a Poisson approximation, that is, it converges to the limiting shape of the classical {\em Richardson model} as $n$ tends to infinity. 

We further investigate how the asymptotic speed of the infection in FCP depends on the variability of the underlying contact-time distribution. For this, we compare the following three models. Two models with precisely one, uniformly distributed, contact time per day and edge, which is either fixed throughout the whole process or independently resampled for each day. The third model is the pure-growth model based on an intensity-one Poisson point process. Importantly, all models are calibrated to have one contact per day in expectation. It turns out that {\em rigidity accelerates infection spread}: The most rigid model with one fixed meeting time per day and edge has the largest limiting set of infected vertices, while the Poisson process based model has the smallest such set, with the resampled meeting time model lying in between the two. In a nutshell, this is related to the well-known waiting-time paradox. 
Hence, a population should organize its contact times as random as possible in order to slow down an infection. Let us mention that, for comparison, the speed of the infection in two-dimensional FPP is increasing {\em if the expected transition time is fixed} but the underlying i.i.d.\ passage-time distribution is made more variable, in the sense of convex ordering, see~\cite{KestenBergFPP,marchand2002strict}. 
In other words, roughly for FPP, {\em more randomness speeds up the infection}, since, in the presence of multiple possible paths, the larger variance provides also particularly fast paths. 

The remainder of the paper is organized as follows. In Section~\ref{sec:result} we present our setup and state our main results about shape theorems for stationary and periodic contact distributions as well as comparisons and connections between FPP and FCP. We conclude with a discussion on further research directions. Section~\ref{sec:proofs} contains all proofs.

\section{Setting and main results} \label{sec:result}
	We consider the lattice \(\Z^d\), $d\ge 1$, and identify it with the {\em nearest-neighbor graph} in the usual way. That is, the vertices are given by the lattice sites and each lattice site is connected to its \(2d\) nearest neighbors. Let \(\{\X_e \subset \R \colon e=\{x,y\}\subset \Z^d\times\Z^d\}\) be a collection of i.i.d.~{\em random closed sets} in \(\R\) indexed by the lattice edges. 
    {\BBB More precisely, considering the Fell-topology on $\mathcal{C}$, the set of closed subsets in $\R$, a random closed set $\X$ is a random variable taking values in $\mathcal{C}$, which is measurable with respect to the Borel-$\sigma$-algebra based on the Fell-topology. This is equivalent to saying that the events
    \begin{equation}\label{eq:Fell-measurability}
        \{ \X \cap K \neq \emptyset \} 
    \end{equation}
    are measurable for all compact sets $K\subset\R$, see e.g.~\cite{molchanov2005theory}.}
    The interpretation {\BBB of the random closed sets} is as follows. Each site represents an agent of a network where two agents occasionally meet or contact each other if they share an edge. The times and durations of these contacts are modeled by the random closed sets assigned to the edges. Canonical examples are Poisson point processes or, more generally, renewal processes and Boolean models thereof. We are interested in the speed of an infection spreading through the network, where the infection is passed from an infected vertex to an uninfected one at the next contact time. Initially, only the agent at the origin \(o\) is assumed to be infected. The time it takes the infection to be passed from a vertex \(x\), infected at time \(t_0\), to some vertex \(y\) is then given by
\begin{equation}\label{eq:definition-D}
		\begin{aligned}
  			D_{t_0}(x,y):= \inf\big\{ 
  				& 
  				t\geq0 \colon \exists \text{ path } \gamma = (x=x_0,x_1,x_2,\dots,x_{n-1},x_n = y) \text{ and }
  				\\& 
  				t_0 \leq t_1 \leq t_2 \leq\dots \leq t_n = t_0 + t \text{ such that }  t_i\in\X_{\{x_{i-1},x_{i}\}} \text{ for all } \, 1\leq i\leq n\big\}.
		\end{aligned}		
	\end{equation} 
	Put differently, one may only traverse neighboring sites along increasing sequences of contact times. Correspondingly, we define the set of reachable vertices from some starting vertex \(x\) during the time interval \([t_0,t_0+t]\) as 
	\begin{equation*}
		I_{t_0}(x,t):=\{y\in \Z^d\colon D_{t_0}(x,y)\leq t\}.
	\end{equation*} 
    {\BBB We discuss measurability of $D_{t_0}(x,y)$, which immediately implies measurability of \(I_{t_0}\), at the beginning of Section \ref{sec:proofs}.}
	As we only assume the origin \(o\) to be infected at initial time \(0\), we abbreviate
	\begin{equation*}
		 D(x,y):=D_0(x,y) \quad \text{ as well as } \quad I(t):=I_0(o,t)=\{y\in\Z^d \colon D(o,y)\leq t\}.
	\end{equation*}
	We call this model {\em first contact percolation} (FCP) on \(\Z^d\). 
    
    \subsection{Stationary contact times}\label{sec:stationaryFCP} 
    For comparison, we recall the basics of {\em first passage percolation} (FPP). Given a probability measure $\mu$ on $[0,\infty]$, we assign to each edge \(e=\{x,y\}\subset\Z^d\) an independent random variable \(X_e\) distributed according to $\mu$. In this case, the transition time between two vertices is defined as the shortest path from \(x\) to \(y\), i.e.,
		\[
			D^{(\mu)}(x,y) := \min \Big\{\sum_{e\in\gamma} X_e\colon \gamma\text{ is a path connecting }x\text{ and }y\Big\}.
		\]  
The set of vertices reached up to time \(t\) from the origin is then given by
		\begin{equation}\label{eq:FPP-infection}
			J^{(\mu)}(t) = \{x\in\Z^d\colon D^{(\mu)}(o,x)\leq t\}.
		\end{equation}
The most prominent example of an FPP model is the {\em Richardson model}, where the $X_e$ are i.i.d.\ exponentially distributed random variables. For an overview and extensive discussion on FPP, we refer the reader to the book~\cite{auffinger2017fpp50years}.
      
Let us note that while \(D^{(\mu)}\) defines a metric, 
this is {\em not true for} \(D_{t_0}\) (as used in FCP). Particularly, the latter is not symmetric as the time it takes to pass the infection from \(x\) to \(y\) may differ considerably from the time it would take the infection to pass from \(y\) to \(x\). Still, \(D_{t_0}\) fulfills a version of the triangle inequality as
\begin{equation}\label{eq:triangle_ineq}
		D_{t_0}(x,z)\leq D_{t_0}(x,y) + D_{D_{t_0}(x,y)}(y,z).
\end{equation}
Another crucial difference between FCP and FPP arises when processes are considered, in which edges are crossed multiple times, for instance a contact process on the graph that includes recoveries, see also the discussion in Section~\ref{sec:discussion}. In this work, we focus on the set of reachable vertices \(I(t)\). In this regard, the link between FCP and FPP is more immanent: {\BBB If the contact times $\X_e$ are stationary under shifts in $\R$, then} $I(t)$ corresponds to an FPP model as each edge is only traversed once.

\begin{proposition}[FCP with stationary contacts is FPP] \label{prop:shape-stationary}
Consider the FCP based on some random closed set $\X$ {\BBB that is stationary under shifts in $\R$}. Then, the process \(\{I(t)\colon t\geq 0\}\) of reachable vertices in the FCP model has the same distribution as \(\{J^{(\mu)}(t)\colon t\geq 0\}\), where
\begin{equation}\label{eq:FCP-FPP-link}
    \mu([0,s]) = \P(\X\cap [0,s]\neq\emptyset)\qquad\text{ for all }s\ge 0.
\end{equation}
\end{proposition}

The proof (and all other proofs) is presented in Section~\ref{sec:proofs}. By now, it is a standard result that the continuum version \(\widetilde{J}^{(\mu)}\) of \(J^{(\mu)}\) satisfies a shape theorem (of the form~\eqref{eq:shape-theorem} below) with limiting shape \(\mathcal{B}^{(\mu)}\) under certain conditions on $\mu$, see e.g.~\cite{auffinger2017fpp50years}. By Proposition~\ref{prop:shape-stationary}, under those conditions, $I$ satisfies the same shape theorem. 
Before we discuss the relation between $\X$ and $\mu$ any further, let us give some examples. 
\begin{example}\label{rem:examples-FCP-FPP} Let us mention some examples for which $\mu$ can be described explicitly.
\begin{enumerate}[(i)]
    \item {\em Poisson point processes}: For $\X$ being the Poisson point process with intensity $\lambda>0$, the corresponding FPP model is the Richardson model based on exponential transition times with parameter $\lambda$.
        \item {\em Renewal processes:} If $\X$ is a stationary simple renewal process with inter-arrival distribution $\nu$ (with finite expectation), then $\mu$ is the distribution of the {\em forward recurrence time}  of $\nu$, i.e., 
        \begin{equation}\label{eq:FPP-from-renewal}
            \mu([0,s])=\int_0^s\nu((x,\infty))\d x\Big/\int_0^\infty\nu((x,\infty))\d x.
        \end{equation}
        In particular, $\mu$ has the non-increasing density $f(s)=\nu((s,\infty))/\int_0^\infty\nu((x,\infty))\d x$.
     \item {\em Shifted lattices}: Consider $\X$ being given by $L(\Z + U)$, $L>0$, where $U$ is a uniform random variable on $[0,1]$. Then, $\mu$ is the uniform distribution on $[0,L]$. 
    \item {\em Boolean models}: Let $\X'$ be a stationary random closed set with associated measure $\mu'$, given by Equation~\eqref{eq:FCP-FPP-link}. 
    Now, let $\X$ be the Boolean model based on $\X'$ with parameter $r>0$, i.e., $\X:=\bigcup_{x\in\X'}[x-r,x+r]$. Then, the corresponding FPP model has a transition-time distribution given by 
        \begin{equation}\label{eq:FPP-from-Boolean-model}
            \mu([0,s]) = \mu'([0,s+2r]).
        \end{equation}
        Indeed, this can be seen after assuming that the balls of radius $r$ are attached ``to the left'' of $\X$, which does not change the distribution due to stationarity.
    \item {\em Cox processes}: Consider the random closed set \(\X\) where we either have \(\X=\R\)  with probability \(1-\e^{-1}\) or \(\X\) is sampled from a Poisson point process of unit intensity otherwise. The corresponding FPP has the transmission-time distribution
        \[
            \mu([0,s])=(1-\e^{-1}) + \e^{-1} (1-\e^{-s}) = 1 - \e^{-(s+1)}.
        \]
        This also corresponds to the FCP given by the Boolean model of a Poisson point process of intensity one and radius \(1/2\). In particular, this shows that different FCPs can correspond to the same FPP, with respect to the process of reachable vertices. 
    \end{enumerate}
\end{example}
In view of the examples just described, our next result describes general properties of $\mu$ defined via~\eqref{eq:FCP-FPP-link} for given $\X$. 

\begin{lemma}[$\mu$ is concave and absolutely continuous]\label{lem:mu-concave}
Let $\X$ be a stationary random closed set and $\mu$ its associated measure as given by Equation \eqref{eq:FCP-FPP-link}. Then, $s\mapsto\mu((0,s])$ is concave. In particular, $\mu$ restricted to $(0,\infty)$ has a Lebesgue density $f$ on $(0,\infty)$, which can be chosen to be non-increasing and right-continuous.
\end{lemma}

What can be said about existence, uniqueness and other properties of stationary random closed sets giving rise to a prescribed transition-time measure $\mu$? In view of Example~\ref{rem:examples-FCP-FPP} Part $(v)$, we see that such $\X$ are not necessarily unique. Moreover, Lemma~\ref{lem:mu-concave} requires concavity of $\mu$. The following result establishes that concavity is also sufficient. Furthermore, we present a construction of some $\X$ via randomly superposed translation-invariant lattices as well as renewal processes.

\begin{proposition}[FCP from FPP]\label{prop:FCPfromFPP}
    Let \(\mu\) be a distribution on \([0,\infty]\) whose restriction to $(0,\infty)$ has a non-increasing, right-continuous Lebesgue density $f$.
    Then, we have the following: 
    \begin{enumerate}[(i)]
        \item There exists a stationary point process \(\X\) such that the FCP under \(\X\) corresponds to the FPP under \(\mu\). Furthermore, $\X$ can be chosen such that each realization is a  shifted lattice of the form \(L(\Z+U)\) for some $L\in[0,\infty]$ and $U\in [0,1]$. Here, we make the convention that \(0\cdot\Z:=\R\) and \(\infty\cdot\Z:=\emptyset\).
        \item If \(\mu\) has no atoms in \(\{0,\infty\}\) and if $f$ is bounded, then there exists a unique stationary and ergodic renewal process \(\X\) such that Equation~\eqref{eq:FCP-FPP-link} is satisfied. 
        \item If \(\mu\) has no atom in \(\infty\), then \(\X\) can be chosen as a Boolean model of a stationary and ergodic renewal process.
    \end{enumerate}
\end{proposition}
Let us note that the stationary sets $\X$ in Part $(i)$ above are not necessarily ergodic. Also, we note that $\X$ in Part $(iii)$ is not unique even within the set of Boolean models based on stationary renewal processes.

So far, we have only discussed processes $\X$ that are stationary with respect to shifts in $\R$, however those processes are not well suited to model the rather periodic behavior of a population with periods of sleep. 

 \subsection{Periodic contact times}\label{sec:periodicFCP}
\begin{figure}[t]
\includegraphics[width=0.33\columnwidth]{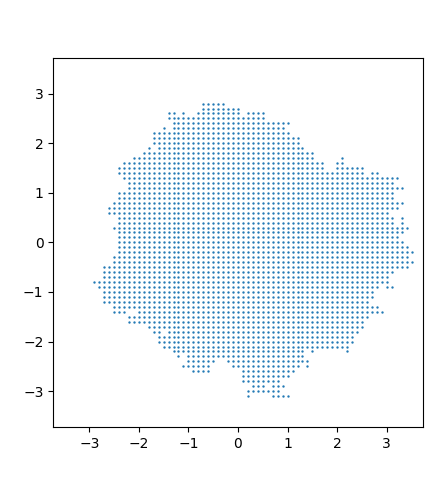}
\includegraphics[width=0.33\columnwidth]{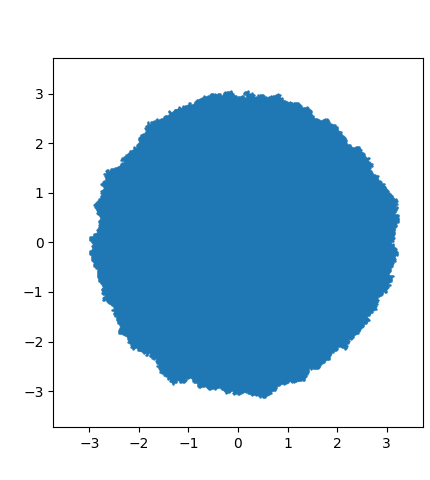}
\includegraphics[width=0.33\columnwidth]{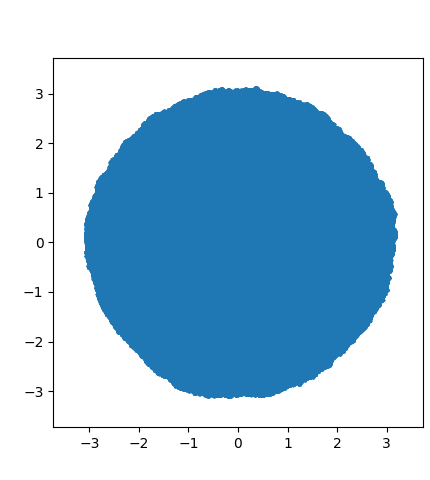}
\caption{Realizations of $t^{-1}I(t)$ for $n=1$ and $t=10,50,190$. One can see the sets slowly converging to a deterministic limit set $\mathcal{B}_1$.}\label{fig:limit-shapes}
\end{figure}

 In order to accommodate periodic behavior, we consider the process based on the following {\em perturbed lattice}  
\begin{equation*}\tag{PL}\label{eq:model_n_Day}
			\X^{(n)} = \bigcup_{k \in \Z} \bigcup_{i=1}^n \{k+U^{(k,i)}\},
\end{equation*}
where the $U^{(k,i)}$ are {\BBB atom-free} i.i.d.\ random variables on $[0,1)$. In other words, $\X^{(n)}$ models a situation in which every interval $[k,k+1)$ contains precisely $n$ i.i.d.\ contacts. {\BBB Multiple contact times cannot coincide since they have no atoms.} Now, every edge is equipped with an independent copy of $\X^{(n)}$ and we write $\{I_{n}(t)\colon t\ge 0\}$ for its process of reachable vertices. Additionally, we write \(\widetilde{I_{n}}(t)\) for its continuum version obtained from \(I_{n}(t)\) by identifying each lattice site \(x\) with the unit box centered around \(x\). 
The following result establishes a shape theorem for \(\widetilde{I_{n}}(t)\) in which the limiting shape does not depend on the distribution of the $U^{(k,i)}$. We present some simulated realizations on $\Z^2$ in Figure~\ref{fig:limit-shapes}.

\begin{thm}[Shape theorem]\label{thm:shape}
		Consider the model based on~\eqref{eq:model_n_Day}. For all \(n\in\N\), there exists a convex, compact {\BBB and non-empty} set $\mathcal{B}_n\subset \R^d$ such that, for all {\BBB atom-free distributions} determining the contact times and $\eps>0$, we have 
		\begin{equation}\label{eq:shape-theorem}
    		\P\big(\exists\  T<\infty\colon (1-\eps)\mathcal{B}_n\subset t^{-1}\widetilde{I_n}(t)\subset(1+\eps)\mathcal{B}_n \text{ for all }t\geq T\big)=1.
		\end{equation}
\end{thm}

The next statement establishes a Poisson-approximation result for the limiting shape. For this, let \(J_R(t)\) denote the set of infected vertices at time $t$ in the 
Richardson model, i.e., the FPP based on $\mu$ being the exponential distribution with parameter one, and \(\Bcal_{\text{R}}\) the limiting shape in the associated shape theorem.  

\begin{thm}[Poisson approximation] \label{thm:approx}
	Let \(I_n(t)\) denote the set of infected vertices at time $t$ and $\Bcal_n$ the limit shape of the corresponding model. Then, there exists a coupling such that almost surely
 \begin{equation*}
     I_n(t/n)\supset J_{\text{R}}(t).
 \end{equation*}
In particular $\Bcal_n/n\supset \Bcal_{\text{R}}$ for every $n\in\N$. Furthermore, with respect to the Hausdorff metric,
		\[
			\Bcal_n/n \to\Bcal_\text{R}\qquad \text{as \,} n\to\infty.
		\]
\end{thm}

Finally, we present a case study that shows that additional randomness in the process $\X$ decreases the set of infected vertices. 

\subsection{Comparison of infection speeds}\label{sec:speedsFCP}
Recall the model $\X=\Z+U$, with $U$ being uniformly distributed in $[0,1)$ from Example~\ref{rem:examples-FCP-FPP} Part~$(iii)$. In this situation, neighboring vertices decide only once independently for a contact time on the first day and then keep that time for all subsequent days. In order to compare this process with the process $\X^{(n)}$ defined in~\eqref{eq:model_n_Day}, let us also consider the multi-contact version of $\X=\Z+U$, given by the {\em stationarized lattice}
\begin{equation*}\tag{SL}\label{eq:modelRigid}
        \X'^{(n)} = \bigcup_{k\in \Z} \bigcup_{i=1}^n \{k+U^{(i)}\} = \bigcup_{i=1}^n \big( \Z + U^{(i)} \big),
    \end{equation*}
where the $U^{(i)}$ are i.i.d.\ uniforms in $[0,1)$. We equip every edge independently with a copy of $\X'^{(n)}$ and write $I'_n(t)$ for the associated process of reachable vertices. 
Clearly, $\X'^{(n)}$ features less randomness than $\X^{(n)}$, which leads to more infections, as shown by the following result.
\begin{proposition}[Randomness slows down]\label{prop:speed-up-through-rigidity}
    Consider, for fixed \(n\in\N\), the model based on~\eqref{eq:modelRigid}. Then, there exists a coupling such that almost surely $I'_n(t) \supset I_n(t)$.
\end{proposition}

Combining Theorem~\ref{thm:approx} and Proposition~\ref{prop:speed-up-through-rigidity}, we have the following ordering of infected sets,
    \begin{equation*}
I'_n(t/n) \supset I_n(t/n)\supset J_R(t)
    \end{equation*}
and thus also the inclusion 
    \begin{equation*}
\Bcal_R\subset\Bcal_n/n\subset \Bcal'_n/n,
    \end{equation*}
{\BBB where the existence of $\Bcal'_n$ follows from Proposition~\ref{prop:shape-stationary} since $\X'^{(n)}$ is stationary and well behaved enough to satisfy a FPP shape theorem, see e.g.~\cite[Theorem~2.16]{auffinger2017fpp50years}.}   Meanwhile, setting $n=1$, all three models have, in expectation, precisely one contact in the interval $[0,1)$ and thus, we can conclude that in our situation, randomness slows down the infection process. As an illustration, let us note that one can explicitly calculate the asymptotic speeds in dimension one for \(n=1\). In this case, one obtains the values $2$ for the FCP model associated to $\X'^{(1)}$,
$\mathrm{e}-1$ for the FCP model associated to $\X^{(1)}$, and $1$ for the Richardson model, underlining the observed ordering.

We conclude this section by mentioning that, as the shape theorem for FPP holds for $I'(t)$, Proposition~\ref{prop:speed-up-through-rigidity} also provides an alternative proof for the finiteness of the speed in the model based on~\eqref{eq:model_n_Day}, see Lemma~\ref{lem:finite-speed-general-n} below.

\subsection{Further discussion and outlook} \label{sec:discussion}
Let us elaborate on some open questions and connections of FCP to other processes. 

\subsubsection{Correspondences and differences between FCP and FPP}
We have seen in Proposition~\ref{prop:FCPfromFPP} that we can find one-to-one correspondences between transition-time distributions $\mu$ in FPP and certain classes of stationary random closed sets $\X$. In that respect, we have considered super-positioned lattices $L(\Z+U)$, renewal processes and Boolean models thereof. However, it should be interesting to see what other properties on $\mu$ are necessary and sufficient to establish links to other classes of stationary random closed sets like, for example, Gibbs point processes or perturbed lattices, in particular ergodic ones. Let us mention that our model~\eqref{eq:model_n_Day} is indeed a perturbed lattice albeit not stationarized. On the other hand, the model~\eqref{eq:modelRigid} is a shifted lattice and fits well into the framework of stationarized renewal processes.

Moreover, let us recall that many different stationary FCP models can lead to the same FPP exploration process. The main difference between FCP and FPP becomes apparent when we consider any characteristic that involves a revisiting of vertices. This begs the question: What other quantities and processes can we look at to easily distinguish the various models? A candidate is introducing recovery, i.e.~considering contact processes. 

\subsubsection{Speed and rigidity}
Theorem~\ref{thm:approx} and Proposition~\ref{prop:speed-up-through-rigidity} have given us a nice ordering of the Richardson model and the FCP models given by $\X$ as in \eqref{eq:model_n_Day} and $\X'$ in \eqref{eq:modelRigid} in terms of their time constants. We have observed that this also matches their ``randomness''. In what generality does this hold and can we make this statement rigorous?

As a small sanity check, let us briefly return to the one-dimensional case $d=1$. Consider the stationarized version of $\X$, that is, $\X+U$ for a uniform random variable on $(0,1)$. In a sense, this model is more random than $\X$, so we would expect its speed to be lower. Indeed, using \eqref{eq:FCP-FPP-link}, we see that its FPP transition-time distribution $\mu$ is given by
\[
    \mu([s,\infty))=(1-s)+s^3 /6 \qquad\text{and}\qquad \mu([1+s,\infty))=(1-s)^3/6
\]
for $s\in[0,1]$. Therefore, its inverse speed is given by $\int_0^\infty \mu([s,\infty))\d s=7/12$. This underlines our assumption on the speeds since $12/7 < \mathrm{e}-1$.

\subsubsection{Beyond stationary and periodic processes $\X$}
Our main results only deal with contact-time processes $\X$ that are either stationary with respect to $\R$ or with respect to some lattice $L\Z$, $L>0$. This has the advantage that the waiting time until the next crossing opportunity of an edge does not depend on the time at which we arrive at an endpoint of that edge. This is no longer the case if we abandon stationarity. Particularly, a shape theorem may no longer hold true, if $\X$, for example, is a renewal process with heavy-tailed inter-arrival times; a case mentioned in~\cite[Proposition 1.3]{hilario2022results}, in the context of contact processes based on renewal processes replacing the usual Poisson processes in the graphical representation. Let $\X$ be a renewal processes started at time $0$ such that, for some $0<\varepsilon<1$, 
\begin{equation*}
\P(\X\cap[t,t+t^\varepsilon]\neq\emptyset)\le t^{-\varepsilon}, \qquad\text{ for all }t\ge t_0,
\end{equation*}
whose existence is guaranteed by~\cite[Proposition~7]{fontes2019contact}. In this case $S_t:=\sup\{\vert y\vert\colon y\in I(t)\}$ behaves sublinearly in the sense that, almost surely for all $0<c<\varepsilon$, 
\begin{equation*}
\limsup_{t\uparrow\infty}t^{c-1}S_t=0. 
\end{equation*}
However, in this example, no stationary version of $\X$ exists. If the inter-arrival times have a finite expectation, then the renewal process will eventually equilibrate and we expect a shape theorem to still hold. 

\subsubsection{Further directions of research}
Natural further directions of research include the 
modeling and analysis of other interacting particle systems via graphical representations not based on Poisson processes, see e.g.~\cite{khoromskaia2014dynamics,concannon2014spatiotemporally} in the context of the exclusion process, \cite{fontes2021contact,fontes2019contact,fontes2020contact,fontes2023renewal,linker2020contact,hilario2022results,seiler2023contact,cardona2024contact} in the context of the contact process, or~\cite{ganguly2022hydrodynamic} for hydrodynamic limits of more general processes.  
In view of the related FPP, several properties can be investigated with respect to FCP, such as
\begin{itemize}
    \item conditions under which we have continuity or monotonicity of the limiting shape with respect to varying random closed set distributions, see~\cite{cox1981continuityFPP,marchand2002strict}, 
    \item large-deviation type results for the probability to reach far-away points atypically fast or slow, see~\cite{chow2003large,basu2017upper},
    \item further properties of the geodesics in FCP such as their fluctuation around the straight line, see~\cite{ahlberg2023chaos,dembin2024superconcentration} or the number of geodesics, see~\cite{hoffman2008geodesics,10.1214/105051604000000503},
    \item properties of the limiting shape for FPP models derived from FCP models, or of periodic FCP models, see~\cite{garet2004asymptotic,f09f98d3-b1ad-3df1-8176-0ddc42f52eab},
\end{itemize}
and many more.   

\section{Proofs} \label{sec:proofs}	
This section is devoted to the proofs of our main results regarding the model based on~\eqref{eq:model_n_Day} but also the  correspondences between FCP and FPP. We start by proving the Propositions~\ref{prop:shape-stationary} and~\ref{prop:FCPfromFPP} in Section~\ref{subsec:correspondences-FCP-FPP} as well as the auxiliary Lemma~\ref{lem:mu-concave}. Together, they provide necessary and sufficient requirements to compare  FPP and certain FCP models.
We prove the shape theorem (\Cref{thm:shape}) in \Cref{sec:proof_shape}. This follows mostly from established arguments. In \Cref{sec:proof_couple}, we prove the monotone coupling between monotonicity the models based on~\eqref{eq:model_n_Day} and~\eqref{eq:modelRigid} (Proposition~\ref{prop:speed-up-through-rigidity}). Additionally, we show convergence of the model based on~\eqref{eq:model_n_Day} to the Richardson model (\Cref{thm:approx}) as $n$ tends to infinity. {\BBB Let us however first verify measurability of the transition times.}

\begin{proof}[Proof of measurability of $D_{t_0}(x,y)$]
    {\BBB
    Given $a\in\R_{\geq0}$, we have
    \[
    	\left\{D_{t_0}(x,y)\leq a\right\} = \bigcup_{\ell\in \N_0}\bigcup_{\substack{\gamma=(x_0,x_1,\dots,x_\ell) \\ x_0=x, x_\ell=y}} \left\{\exists t_0\leq t_1 \leq \dots \leq t_\ell \leq t_0+a \colon t_i\in \X_{(x_{i-1},x_i)} \,\forall 1\leq i \leq \ell\right\}.
    \]
    It therefore suffices to study measurability along the paths $\gamma$. We argue by induction in $\ell$ that, for every path $\gamma=(x_0,x_1,\dots,x_\ell)$, every $a\in\R_{\geq0}$ and every $t_0\in\R$, all events of the form
    \[
    	A_{t_0}(\gamma,a):=\left\{\exists t_0\leq t_1 \leq \dots \leq t_\ell \leq t_0+a \colon t_i\in \X_{(x_{i-1},x_i)} \,\forall 1\leq i \leq \ell\right\} 
    \] 
    are measurable. This is indeed the case for $\ell=1$, as
    \[
    	A_{t_0}(\gamma,a) = \left\{\exists t_0\leq t_1 \leq t_0+a \colon t_1\in \X_{(x_0,x_1)}\right\} = \{\X_{(x_0,x_1)}\cap [t_0,t_0+a]\neq\emptyset\},
    \]
    which is measurable by the Characterization~\eqref{eq:Fell-measurability}. Consider now a path $\gamma=(x_0,x_1,\dots,x_\ell,x_{\ell+1})$. Writing $\gamma_-=(x_0,x_1,\dots,x_\ell)$, we have
    \begin{equation*}
        \begin{aligned}
            A_{t_0}(\gamma,a) 
            &= \bigcup_{b\in[0,a]} \left( A_{t_0}(\gamma_-,b)\cap \{\X_{x_\ell,x_{\ell+1}} \cap [t_0+b, t_0+a]\} \right)\\
            &= \bigcap_{k\ge 1} \bigcup_{b\in[0,a]\cap\mathbb{Q}} \left( A_{t_0}(\gamma_-,b+\tfrac1k)\cap \{\X_{x_\ell,x_{\ell+1}} \cap [t_0+b-\tfrac1k, t_0+a]\} \right),
        \end{aligned}
    \end{equation*}
    which is measurable by the induction hypothesis and again \eqref{eq:Fell-measurability}. Therefore, $D_{t_0}(x,y)$ is measurable as a countable intersection over a countable union of measurable events.
    }
\end{proof}

\subsection{Proofs of FCP-FPP correspondences}\label{subsec:correspondences-FCP-FPP}

\begin{proof}[Proof of Proposition~\ref{prop:shape-stationary}]
	Observe that $\{I(t)\colon t\geq 0\}$ can be viewed as an exploration process of \(\Z^d\), starting from the origin, that we now describe.
    \BBB{
    	At time zero, we start with the origin, mark it as explored and set \(I(0)=\{o\}\). We reveal the origin's neighboring edges in the sense that we reveal their respective waiting times until the next contact, which are independent and distributed according to \(\mu\), defined in~\eqref{eq:FCP-FPP-link}. Let us denote by \(\tau_1^{0},\tau_2^{0},\dots\) these next contact times, ordered from smallest value to largest, with an appropriate tie breaker if necessary. At time \(\tau^{0}_1\), we mark the yet unexplored end vertex of the corresponding edge, say \(v_1\), as being explored and set \(I(\tau^{0}_1)=I(\tau^{0}_1-)\cup \{v_1\}=\{o,v_1\}\). We reveal the boundary edges of \(v_1\), which are all edges connecting \(v_1\) to a yet unrevealed vertex and reveal the waiting time until the next contact. These newly revealed waiting times are independent and again distributed according to~\(\mu\). We update the previous waiting times \(\tau_2^{0}-\tau_1^{0}, \tau_3^{0}-\tau_1^{0},\dots\) and add the newly revealed waiting times, and sort all of them from smallest to largest value to obtain \(\tau_1^{1},\tau_2^{1},\dots\). We proceed in the obvious way with the extension that, at later stages, we may explore a vertex that has edges incident to vertices that have already been explored at an earlier time. If this is the case, we remove the corresponding waiting time from the (ordered) list of waiting times. 
    
    	Let us consider now the exploration process \(\{J^{(\mu)}(t)\colon t\geq 0\}\) corresponding to the FPP model with i.i.d.\ transition times with distribution \(\mu\). It is straightforward to describe \(J\) as an exploration of \(\Z^d\) just like before. Again we start with \(J^{(\mu)}(0)=\{o\}\) and reveal the transition times of the incident edges, order them from smallest value to largest and denote them by \(\sigma_1^0,\sigma_2^0,\dots\). Then, at time \(\sigma_1^0\), we explore the corresponding neighbor and set \(J^{(\mu)}(\sigma_1^0)=\{o,v_1\}\). We further update the transition times \(\sigma_2^0-\sigma_1^0,\dots\), reveal the new transition times, incident to \(v_1\), sort them from smallest value to largest to obtain \(\sigma_1^1,\sigma_2^1,\dots\). Again, we proceed in the obvious way. 
    	
    	As the waiting times \(\tau_i^j\) and \(\sigma_i^j\) are clearly identically distributed, these exploration describes a coupling between \(I\) and \(J^{(\mu)}\) such that \(I(t)=J^{(\mu)}(t)\) for each \(t\geq 0\) that ultimately proves Proposition~\ref{prop:shape-stationary}.
    }
\end{proof}

\begin{proof}[Proof of Lemma~\ref{lem:mu-concave}]
    The proof is identical to the one in \cite[Lemma~3.1.1]{beutler1966theory} for stationary point processes, which we recapitulate here for convenience. Let $h>0$ and $0<a<b<\infty$. We first observe that
    \[
    \begin{aligned}
        \mu((a,a+h]) 
        &= \P(\X\cap[0,a]=\emptyset\text{ and }\X\cap(a,a+h]\neq\emptyset)\\
        &= \P(\X\cap[b-a,b]=\emptyset\text{ and }\X\cap(b,b+h]\neq\emptyset)\\
        &\geq \P(\X\cap[0,b]=\emptyset\text{ and }\X\cap(b,b+h]\neq\emptyset) = \mu((b,b+h]),
    \end{aligned}
    \]
    where the second equality follows from stationarity of $\X$. Given arbitrary $0<s<t<\infty$ and choosing $a=s,\,b=(s+t)/2$ and $h=(t-s)/2$ yields
    \[
        \mu\big((s,(s+t)/2]\big)\geq \mu\big(((s+t)/2,t]\big)
        \iff \mu\big((0,(s+t)/2]\big)\geq 2^{-1} \big( \mu((0,s])+\mu((0,t])\big),
    \]
    which shows concavity of the function $F\colon(0,\infty)\to[0,1],\,s\mapsto \mu((0,s])$. Hence, a right derivative $f(s)$ exists for all $s>0$. Additionally, $f$ is decreasing and $F(c)-F(b) \leq f(a)|c-b|$ for all $0<a<b<c$. Thus, $F$ is absolutely continuous on any interval $[\delta,\infty),\,\delta>0$ and hence
    \[
        F(b) = F(a) + \int_a^b f(s)\d s,
    \]
    concluding the proof by taking the limit \(a\downarrow 0\). 
\end{proof}

The following observation will be used in the proof of~Proposition~\ref{prop:FCPfromFPP}.
\begin{lemma}[Limit behavior of $f$]\label{lem:f-at-0-infty}
    Let $f\colon (0,\infty)\to(0,\infty)$ be monotone with $\int_0^\infty f(x)\d x<\infty$. Then,
    \[
    \lim_{x\downarrow 0}xf(x)=\lim_{x\uparrow\infty}xf(x)=0.
    \]
\end{lemma}
\begin{proof}
    It follows immediately that $f$ must be non-increasing.
    Next, assume $\limsup_{x\uparrow\infty}xf(x)=:a>0$. Then, there exists an increasing sequence $(x_k)_{k\in\N},\,x_k\uparrow\infty$ such that $x_kf(x_k)\geq a$, i.e.,~$f(x_k)\geq a/x_k$. By choosing a subsequence, we may assume, without loss of generality, that $x_{k-1}/x_k<1/2$. Then, due to monotonicity of $f$, we have that
    $$\int_0^\infty f(x)\d x\geq \sum_{k=2}^\infty (x_k-x_{k-1})f(x_k)
    \geq \sum_{k=2}^\infty (1-x_{k-1}/x_k)a
    \geq \sum_{k=2}^\infty a/2 = \infty,$$
    which shows that $\limsup_{x\uparrow\infty}xf(x)=0$ whenever $\int_0^\infty f(x)\d x<\infty$. Showing $\limsup_{x\downarrow 0}xf(x)=0$ works analogously,
    which finishes the proof.
\end{proof}

\begin{proof}[Proof of Proposition~\ref{prop:FCPfromFPP}]
We prove the statements individually. 

    {\bf Part $(i)$:} Let us first additionally assume that $\mu(\{0,\infty\})=0$. Take $L$ to be distributed according to a probability distribution $\nu$ on $(0,\infty)$ to be specified later and take $U$ to be a uniform random variable on $(0,1)$ independent of $L$. Our random closed set is going to be $\X:=L(\Z+U)$, which is clearly stationary but not necessarily ergodic. Let us first define $\nu$. 
    Since $f$ is non-increasing on $(0,\infty)$, we may define a measure $Q$ on $(0,\infty)$ via 
    $$Q((a,b]):=f(a)-f(b),\qquad\text{ for all }0<a<b<\infty\,,$$
    {\BBB which we may do since $f$ is right-continuous.}
    Then, we define $\nu$ to be the measure that is absolutely continuous with respect to $Q$ with density $x\mapsto x$, i.e., 
    \[
        \nu(\d x) := x Q(\d x).
    \]
    Let us first check that $\nu$ is indeed a probability measure. For this, we first see that Lemma~\ref{lem:f-at-0-infty} above applies to $f$, i.e.,~$\lim_{x\downarrow0}xf(x)=\lim_{x\uparrow\infty}xf(x)=0$. 
    Further, we obtain by the integration-by-parts formula for Stieltjes integrals for $0<a<b<\infty$, 
    \[
        \int_a^b x \, Q(\d x) = -\int_a^b x \, \d f(x) = af(a)-bf(b)+\int_a^b f(x) \, \d x = a f(a)-b f(b) + \mu([a,b]).
    \]
    Hence,
    \begin{equation*}
        \nu((0,\infty)) = \int_0^\infty x \, Q(\d x) = \lim_{x\downarrow 0}xf(x) - \lim_{x\uparrow\infty}xf(x) + \mu((0,\infty)) = 1.
    \end{equation*}
    Let us verify Equation \eqref{eq:FCP-FPP-link} next. 
    The procedure is as follows: Given $s>0$, the set $\X=L(\Z+U)$ will always have a point in $[0,s)$ if $s\geq L$. Otherwise, we need $LU\leq s$.
    Using integration by parts again, we infer
    \[
    \begin{aligned}
            \P(\X\cap[0,s)\neq\emptyset) 
            &= \P(L(\Z+U)\cap[0,s)\neq\emptyset \text{ and }  L> s) + \P(L(\Z+U)\cap[0,s)\neq\emptyset \text{ and } L\le  s)\\
            &= \int_s^\infty  \P(xU\leq s) \, \nu(\d x) + \int_0^s \, \nu(\d x)= \int_s^\infty x\frac{s}{x} \, Q(\d x) + \int_0^s x \, Q(\d x)\\
            &= sQ((s,\infty)) - \int_0^s x \, \d f(x)= sf(s) - sf(s) + \lim_{x\downarrow 0} xf(x) + \int_0^s f(x) \, \d x = \mu([0,s]).
    \end{aligned} 
    \]
Now, let us consider the general case of $\mu(\{0,\infty\})\geq0$. In this case, we take 
    $$\nu := \mu(\{0\})\delta_0 + \mu(\{\infty\})\delta_\infty + (1-\mu(\{0,\infty\}))\nu' $$
   and interpret $0\cdot\Z=\R$ and $\infty\cdot \Z:=\emptyset$. In the case that $\mu(\{0,\infty\})<1$, we choose $\nu'$ as in the case above for the measure $\mu'([a,b]):=(1-\mu(\{0,\infty\}))^{-1} \mu([a,b]\cap (0,\infty))$, $0\le a<b\le \infty$, where the latter is indeed a probability measure satisfying the conditions required in the above first part. Thus, $\X=L(\Z+U)$ satisfies Equation~\eqref{eq:FCP-FPP-link} since $\mu(\{0\})=\P(0\in\X)=\P(L=0)$ and analogously for $\mu(\{\infty\})$.

    \medskip
    {\bf Part $(ii)$:} Assume $\mu(\{0,\infty\})=0$ and $f(0)<\infty$ for the right-continuous, monotone Lebesgue density of $\mu$. We aim to finding a stationarized (simple) renewal process $\X$ with inter-arrival distribution $\nu$ generating $\mu$. Inverting Equation~\eqref{eq:FPP-from-renewal} yields
    \[
        \nu((s,\infty)) := f(0)^{-1} f(s).
    \]
    The validity of Equation~\eqref{eq:FCP-FPP-link} is easy to check.

    \medskip
    {\bf Part $(iii)$:} We allow for $\mu(\{0\})=:c>0$ and want to use a Boolean model of a stationary and simple renewal process instead. If $c=1$, simply take $\nu := \delta_1$ and radius $r=\infty$. Otherwise, $\mu$ has a Lebesgue density $f$ on $(0,\infty)$. Consider the non-increasing, right-continuous function $f^*\colon [0,\infty)\to[0,\infty)$ as 
    $$f^*(x):=f\big(\max(0,x- c/f(0))\big).$$ 
    We see that $f^*$ is a probability density as 
    $$\int_0^\infty f^*(x)\d x = \int_0^\infty f(x)\d x + \frac{c}{f(0)}f(0) = \mu((0,\infty))+\mu(\{0\})=1.$$
    Choosing $\nu$ as in Part~$(ii)$ but according to $f^*$ and taking $r=c/(2f(0))$ proves the claim by Equation~\eqref{eq:FPP-from-Boolean-model}. Note that many other candidates for $f^*$ exist that yield a smaller radius $r$, hence $\X$ is not unique. 
\end{proof}

\subsection{Proof of shape theorem}\label{sec:proof_shape} 
During the whole section, we fix some \(n\in\N\) for the number of contact times per day and may drop it from the notation whenever convenient. Furthermore, we assume that the contact times per day are i.i.d.\ {\BBB according to some diffusive measure. However, we will see that the particular choice of the measure has no influence on the result.} 

The proof of the shape theorem then follows established arguments, in which the first step is to establish linear speed in a fixed direction. That is, for each fixed direction \(x\in\R^d\), there exists some \(\varphi(x)\) such that the time it takes the infection to reach a site close to \(t x\) is roughly \(t\varphi(x)\). The quantity \(\varphi(x)\) is hence also called the {\em inverse speed} in direction \(x\). For $x\in\R^d$, let us write \([x]\) for its closest neighbor in \(\Z^d\). In case of ambiguity, any tie-breaker will do.

\begin{lemma}[Existence of inverse speed for lattice points]\label{lem:existence-inverse-speed}
	For all $x\in \Z^{d}$, there exists $\varphi(x):=\varphi^{(n)}(x)\in[0,\infty)$ such that almost surely
	\[
		\lim_{t\uparrow\infty}t^{-1}D(o,[tx])=\varphi(x).
	\]
	Furthermore,
	\begin{equation}\label{eq:characterization-of-phi}
    	\varphi(x)=\lim_{t\uparrow\infty} t^{-1}\E D(o,[tx]) = \inf_{t\geq 1} t^{-1}\E D(o,[tx]).
	\end{equation}
\end{lemma}
\begin{proof} 
In order to make use of the independence between days, it will be useful to consider an integer-valued version of the travel duration $D$, defined as
\begin{equation*}
	\widetilde{D}_{t_0}(y,z):=\lceil D_{t_0}(y,z)\rceil=\inf\{t\in\N\colon z\in I_{t_0}(y,t)\},
\end{equation*}
for a time \(t_0\geq 0\). That is, the infection starting in \(x\) at time \(t_0\) reaches vertex \(y\) in \(\widetilde{D}_{t_0}(y,z)\) days or during the \(\lfloor t_0\rfloor + \widetilde{D}_{t_0}(y,z)\)-th day of the process respectively. Note that the limits $D(0,[tx])/t$ and $\widetilde D_0(0,[tx])/t$ (as $t\uparrow\infty$) coincide if they exist. Given $0\le s\le t \in\N$, we define the random variables
\[
	X_{s,t}:=\widetilde{D}_{\widetilde{D}_0(o,[sx])}([sx],[tx]),
\]
representing the number of days required for the infection to spread from $[sx]$ to $[tx]$ when started the first day after the infection has reached \([sx]\) from the origin. In the remainder of the proof, we shall identify \(sx\) and \(tx\) with $[sx]$ and $[tx]$ respectively to ease notation. The proof is now concluded by applying the subadditive ergodic theorem~\cite[Theorem~2.6]{liggett1985interacting} in the version of~\cite[Theorem~4.2]{deijfen2003asymptotic}. For this, we have to verify the following conditions:
\begin{enumerate}[(i)]
	\item \(X_{0,t}\leq X_{0,s}+X_{s,t}\),
	\item \(0\leq \E X_{0,t}<\infty\) for every \(t\),
	\item the distribution of \(\{X_{s,s+k}\colon k\in\N\}\) does not depend on \(s\), and
	\item \(\limsup_{t\uparrow\infty} X_{0,tk}/t\leq \E X_{0,k}\) for each \(k\in\N\).
\end{enumerate}
Under the assumption that the Conditions (i)--(iv) are satisfied, we infer from~\cite[Theorem~4.2]{deijfen2003asymptotic}
\[
\varphi(x):=\lim_{t\to\infty} t^{-1}\E X_{0,t} = \inf_{t\geq 1} t^{-1}\E X_{0,t},
\]
implying~\eqref{eq:characterization-of-phi}, and \(\varphi(x)=\lim_{t\to\infty} t^{-1} X_{0,t}\) almost surely. It hence remains to verify these conditions.    

\begin{description}
	\item[On (i)] This is a direct consequence of the triangle inequality~\eqref{eq:triangle_ineq}.
	\item[On (ii)] This is immediate from the trivial bounds \(0\leq X_{0,t}\leq t|x|_1\). 
	\item[On (iii)] Note that we have for all measurable sets $A\subset \N^\N$
    \begin{equation*}
	\begin{aligned}
    \P\big(X_{s,s+1},X_{s,s+2},\dots \in A\big)
    &=\P\big(\widetilde D_{\widetilde D_0(o,sx)}(sx,(s+1)x),\widetilde D_{\widetilde D_0(o,sx)}(sx,(s+2)x),\dots \in A\big)
    \\
    &=\sum_{j\ge 0}\P\big(\widetilde D_j(sx,(s+1)x),\widetilde D_j(sx,(s+2)x),\dots \in A,\widetilde D_0(0,sx)=j\big)
    			\\ 
                &=\sum_{j\ge 0}\P\big(\widetilde D_j(sx,(s+1)x),\widetilde D_j(sx,(s+2)x),\dots \in A\big)\P(\widetilde D_0(0,sx)=j)\\ 
                &=\sum_{j\ge 0}\P\big(\widetilde D_0(0,x),\widetilde D_0(0,2x\big),\dots \in A\big)\P(\widetilde D_0(0,sx)=j)\\ 
                & =\P\big(X_{0,1},X_{0,2},\ldots \in A\big),
		\end{aligned}
		\end{equation*}
		where we used the independence of the contact times with respect to different days as well as the lattice shift invariance of the system on every individual day in the penultimate step.

	\item[On (iv)] 
		Define recursively
		\begin{equation*}
			Y^k_1:=\widetilde D_0(0,kx)=X_{0,k}
					\qquad \text{ and }\qquad 
			Y^k_i:=\widetilde D_{Y_1^k+\cdots+Y_{i-1}^k}((i-1)kx,ikx)
		\end{equation*}
		and note that for every realization
		\begin{equation*}
			X_{0,nk}\le \sum_{i=1}^nY^k_i.
		\end{equation*}
		{\BBB Thus, the claim follows by the strong law of large numbers after verifying that the \(Y_i^k\) form an i.i.d.\ sequence. First observe that for $t,a\in\N$, the events 
        \[
        	A_i(t,a):=\{ \widetilde D_{t}((i-1)kx,ikx) = a \} 
        \]
        only depend on the days $t,t+1,\dots,t+a-1$. Also note that the probability of occurrence of $A_i(t,a)$ depends neither on $i$, due to spatial stationarity on $\Z^d$, nor on $t$, due to temporal stationarity on $\Z$. Now, let $r\in\N$ and then $a_1,\dots,a_r\in\N$. It follows that the events $A_i(a_1+\dots+a_{i-1},a_i)$ are independent for different $i$. The event $\{Y_1+\dots+Y_{r-1}=\ell\}$ only depends on the days up to $\ell-1$, and hence, 
        \begin{equation*}
            \begin{aligned}
                \P\left( Y_r^k = a_r     \right) 
                &
                	= \sum_{\ell=1}^\infty \P\left( Y_1+\dots+Y_{r-1}=\ell \text{ and } Y_r^k = a_r      \right)
					= \sum_{\ell=1}^\infty \P\left( Y_1+\dots+Y_{r-1}=\ell \text{ and } A_r(\ell,a_r)     \right)
                \\ &
                	= \sum_{\ell=1}^\infty \P\left( Y_1+\dots+Y_{r-1}=\ell\right)\, \P\left( A_r(\ell,a_r)     \right)
                	= \sum_{\ell=1}^\infty \P\left( Y_1+\dots+Y_{r-1}=\ell\right)\, \P\left( A_1(0,a_r)     \right)
                \\ &
                = \P\left(A_1(0,a_r)\right),
            \end{aligned}
        \end{equation*}
        where we used independence in the second to last step. This implies that the $Y_i^k$ are identically distributed. Their independence follows from        
        \begin{equation*}
            \begin{aligned}
                \P\Big(\bigcap_{i=1}^r \{Y_i^k = a_i\}\Big) 
                &= \P\Big(\bigcap_{i=1}^r \{\widetilde D_{a_1+\dots+a_{i-1}}((i-1)kx,ikx) = a_i\} \Big)
                = \P\Big(\bigcap_{i=1}^r A_i({a_1+\dots+a_{i-1}}, a_i)\Big)
                \\
                &= \prod_{i=1}^r \P\big(A_i({a_1+\dots+a_{i-1}}, a_i)\big)
                = \prod_{i=1}^r \P\left( A_1(0, a_i)     \right)
                = \prod_{i=1}^r \P\left( Y_i^k = a_i     \right),
            \end{aligned}
        \end{equation*}
        finishing the proof of Claim~(iv) and therefore the proof of Lemma~\ref{lem:existence-inverse-speed}.

        }
\end{description}
\end{proof} 

Next, we extend the established inverse speed in \(\Z^d\) directions to {\em all} directions.

\begin{lemma}[Existence of inverse speed for general directions]\label{lem:existence-inverse-speed-general-direction} 
	Under the assumptions and definitions of Lemma \ref{lem:existence-inverse-speed}, almost surely, for all $x\in\R^d$, the limit {\BBB $\varphi(x)=\varphi^{(n)}(x) := \lim_{t\uparrow\infty}t^{-1}D(o,[tx])$} exists. Furthermore, $\varphi\colon\R^d\to[0,\infty)$ is almost surely continuous.
\end{lemma}
\begin{proof}
    Since we travel at least one edge per day, the reverse triangle inequality implies \(|D_0(o,y)-D_0(o,z)|\leq |z-y|_1\)
    for every $y,z\in\Z^d$. In particular, this yields for arbitrary reals  $y,z\in\R^d$, 
    \begin{equation}\label{eq:reverse-triangle-real-points}
    	|D_0(o,[y])-D_0(o,[z])|\leq |z-y|_1 + 2d.    
    \end{equation}
    Let us first consider the case $x\in\mathbb{Q}^d$. Then, there exists some $s\in\N$ such that $sx\in\Z^d$. We will show that $\varphi(x)$ exists and 
    \[
\varphi(x)=\lim_{t\uparrow\infty} t^{-1}D_0(o,[tx])=\lim_{t\uparrow\infty} (st)^{-1}D_0(o,stx)=s^{-1}\varphi(sx).
    \]
    By Lemma \ref{lem:existence-inverse-speed}, \(\varphi(sx)\) exists almost surely. Furthermore, for each \(t>0\), we pick \(k=k(t)\) such that $t\in [(k-1)s,ks)$. Then, using the triangle inequality~\eqref{eq:triangle_ineq} and~\eqref{eq:reverse-triangle-real-points}, we infer
    \[
    	\begin{aligned}
        	\lim_{t\uparrow\infty} &  |s^{-1}  \varphi(sx) - t^{-1}D_0(o,[tx])| 
        	\\ &
        		\leq \lim_{t\uparrow\infty} \big|s^{-1}\varphi(sx) - (ks)^{-1}D_0(o,ksx)\big| + \lim_{t\uparrow\infty} \big|(ks)^{-1}D_0(o,ksx) - t^{-1}D_0(o,[tx])\big|
        	\\ &
        		\leq 0 + \lim_{t\uparrow\infty}(ks)^{-1}\big|D_0(o,ksx) - D_0(o,[tx])\big| + \lim_{t\uparrow\infty}\big|(ks)^{-1} - t^{-1}\big| D_0(o,[tx])
        	\\ &
        		\leq \lim_{t\uparrow\infty} (ks)^{-1} \big[(ks-t)|x|_1 + 2d\big] + \lim_{t\uparrow\infty} \big|(ks)^{-1} - t^{-1}\big|\big(|tx|_1 + d\big) 
        	\\ &
        		\leq \lim_{t\uparrow\infty} (ks)^{-1} \big[s|x|_1 + 2d\big] + \lim_{t\uparrow\infty} (\tfrac{1}{(k-1)s}-\tfrac{1}{ks}) \big(ks|x|_1 + d\big) 
        	\\ &
        	= 0,
    	\end{aligned}
    \]
    using $k(t)\to\infty$ as $t\to\infty$. This shows the claim for each rational direction $x\in\mathbb{Q}^d$. The proof extends to real directions $x\in\R^d$ as
    \[
    	 \lim_{t\uparrow\infty} t^{-1}D_0(o,[tx])
    \]
    is continuous in $x$ since~\eqref{eq:reverse-triangle-real-points} yields
    \[
        t^{-1}|D_0(o,[tx])-D_0(o,[ty])| \leq |y-x|_1 + 2d/t,
    \]
    which concludes the proof.
\end{proof}

Having established the inverse speed in all directions, we proceed by showing its positivity, i.e.,\ the infection has finite speed. Let us mention that, unlike before, the dependence on \(n\) becomes relevant now, which we shall emphasize in the notation. In order to prove finite speed, we require the notion of a path's {\em traveling time}: Given a path $\gamma = (x_0,\dots,x_k)$, we {\BBB write $|\gamma|=k$ for its length and }define its traveling time from $t_0\geq 0$ by 
\begin{equation}\label{eq:def:traveltime}
  T^{(n)}_{t_0}(\gamma):= \inf\big\{t\geq 0\colon  \exists t_0 \leq t_1 \leq\dots \leq t_k = t_0 + t \text{ such that }
   t_i\in\X_{(x_{i-1},x_{i})}^{(n)} \text{ for all } 0<i \leq k\big\}  
\end{equation}
and write $T^{(n)}(\g):=T_0^{(n)}(\g)$ as usual.
Clearly, the minimum in \eqref{eq:def:traveltime} is always realized {\BBB since the $\X^{(n)}_e$ are closed, but possibly taking the value $\infty$}. {\BBB Let us make the following observation for  self-avoiding paths: Each edge $e$ and hence each $\X_e^{(n)}$ is visited at most once. As such, we see from (\ref{eq:def:traveltime}) that the distribution of $T^{(n)}_0(\gamma)$ only depends on $|\gamma|$, but not on the specific path itself.}

\begin{lemma}[Finite speed]\label{lem:finite-speed-general-n}
	For every \(n\in\N\) and $x\in\R^d$ with $|x|_\infty\geq1$, we have $\varphi^{(n)}(x)\geq (2d n \e)^{-1}$. In particular, the speed is finite.
\end{lemma}
\begin{proof}
	We aim to show that for all $a<(2dn\e)^{-1}$, we have $\lim_{t\to\infty}\P(D^{(n)}(o,[tx])\ge at) = 1$. First, a union bound over self-avoiding paths $\gamma$ starting in the origin $o$ yields 
	\begin{equation}
		\P(D^{(n)}(o,[tx])<at)\le \sum_{k\ge t}\sum_{\g\colon |\g|=k}\P(T^{(n)}(\g)<at).
	\end{equation}
In order to bound the probability on the right-hand side, let us introduce the {\BBB infinite, self-avoiding path  \(\gamma_\infty=(o,e_1,2e_1,3e_1,\dots)\).  Since the distribution of $T^{(n)}(\g)$ only depends on $|\gamma|=k$, we have}
\[
	\P(T^{(n)}(\g)<at) \leq \P(\text{At least the first } k  \text{ edges of } \gamma_\infty \text{ are traversed in } \lceil at\rceil \text{ days}). 
\]
Let $\tau^{(n)}_i$ be the (random) number of edges that $\gamma_\infty$ travels during the \(i\)-th day. {\BBB We next show that the \(\tau^{(n)}_i\), \(i=1,2,\dots\), form an i.i.d.\ sequence, using a procedure analogous to the proof of Lemma~\ref{lem:existence-inverse-speed}~(iv). Given $s,k\in\N$, consider the event that the infection, starting in $se_1$ on day $i$ and being passed in a straight line, reaches exactly $(s+k)e_1$ and no further. That is,
    
\[
    A_i(s,k):=\left\{T^{(n)}_i\left((se_1,(s+1)e_1,\dots,(s+k)e_1) \right)\leq 1\right\} \cap \left\{T^{(n)}_i\left((se_1,(s+1)e_1,\dots,(s+k+1)e_1) \right) > 1\right\}.
\]
Observe that $A_i(s,k)$ only depends on the random variables 
$\{U_e^{(i,j)}\colon e\in E,j\in\{1,\dots,n\}\}$
that make up the contact times of the edges 
$e\in E$ on the $i$-th day in \eqref{eq:model_n_Day}.
The independence of the events $A_i(s,k)$, \(i=1,2\dots,\) is then a direct consequence of the independence of the contact times of each day.  Additionally, the probability of occurrence of $A_i(s,k)$ does not depend on $s$ by spatial stationarity. Consider now $k_1,\dots, k_r\in\N$ for some $r\in\N$. It follows that $\sum^{r-1}_{i=1}\tau_i$ depends only on 
$\{U_e^{(j,i))}\colon e\in E,j\in\Z,\,i\leq r-1\}$, which is independent of $A_r(s,k)$ for any $s,k\in\N$. Thus,
\begin{equation*}
    \begin{aligned}
        \P(\tau_r^{(n)}=k_r)
        &=\sum_{l=0}^\infty \P\Big(  \sum^{r-1}_{i=1}\tau_i=l \text{ and } \tau_r=k_r   \Big)
        =\sum_{l=0}^\infty \P\Big(  \sum^{r-1}_{i=1}\tau_i=l \text{ and } A_r(l,k_r) \Big)\\
        &= \sum_{l=0}^\infty \P\Big(  \sum^{r-1}_{i=1}\tau_i=l\Big)\, \P\big(A_r(l,k_r)   \big)
        = \sum_{l=0}^\infty \P\Big(  \sum^{r-1}_{i=1}\tau_i=l\Big)\, \P\big(A_1(0,k_r)   \big)
        = \P\big(A_1(0,k_r)   \big).
    \end{aligned}
\end{equation*}
However, this shows that $(\tau_i)_{i\in\N}$ is an i.i.d.\ sequence, as
\begin{equation*}
    \begin{aligned}
\P\big(\tau^{(n)}_i=k_i\,\,\forall i\leq r\big)
        &= \P\Big( \bigcap_{i=1}^r A_i\Big(\sum^{i-1}_{l=1}k_l,k_i\Big)   \Big)
        =\prod_{i=1}^r\P\Big(  A_i\Big(\sum^{i-1}_{l=1}k_l,k_i\Big)\Big)\\
        &= \prod_{i=1}^r\P\big(  A_1\left(0,k_i\right)\big)
        = \prod_{i=1}^r \P(\tau^{(n)}_i=k_i).
    \end{aligned}
\end{equation*}
}
Further, note that
\[
	\P(\tau^{(n)}_1\geq j) \leq n^j/j! \quad \text{ for all } j\in\N,
\]
as there are \(n\) contact times on each of the \(j\) edges, yielding \(n^j\) possibilities to choose one contact time from each edge with each such choice forming an increasing sequence with probability \(1/j!\). Note that $\tau^{(n)}_1$ has all exponential moments and in particular
\[
	\E\tau^{(n)}_1 \leq \e^n \quad \text{ and } \quad  \E \e^{s\tau^{(n)}_1} \leq  \exp(n\e^s)<\infty \ \text{ for any }s\in\R.
\]
Hence, by an exponential moment bound
\begin{equation}\label{eq:prob-path-too-fast-1}
	\begin{aligned}
    	\P(T^{(n)}(\gamma)\leq at)
    	&
    		\leq \P\Big(\sum_{i=1}^{\lceil at\rceil} \tau^{(n)}_i \ge k\Big)
   			\leq \prod\limits_{i=1}^{\lceil at\rceil}\E \exp\big(s\tau^{(n)}_i\big)\cdot{\exp(-ks)}
    		\leq \exp{(\lceil at\rceil n\e^s-ks)}.
	\end{aligned}
\end{equation}
Note that there are no more than $(2d)^k$ paths of length $k$. Choosing $s=\log(2d)+1$ yields

\begin{equation}\label{eq:prob-path-too-fast-2}
	\begin{aligned}
    	\P(D(0,[tx])<at)& 
    	\leq \sum\limits_{k\geq t}\sum\limits_{\gamma:|\gamma|=k} \exp{(\lceil at\rceil n\e^s-ks)} 
    	\\ &
			\leq\sum\limits_{k\geq t}(2d)^k\exp(-k[\log(2d)+1]+\lceil at\rceil n\e^{\log(2d)+1})
		\\ &
			\leq \sum\limits_{k\geq t}\e^{-k}\exp(2d\e (at+1) n)\leq\e^{-t} \frac{\e^{1+2den}}{\e-1} \exp(2d\e at n),
	\end{aligned}
\end{equation}
which tends to zero as $t\to\infty$ for every $a<1/(2dn\e)$. This concludes the proof.
\end{proof}
\begin{cor}[Exponential decay]\label{cor:exp-decay-general-n} 
	For every $a>0$ and $n,t\in\N$, we have
    {\BBB
    \begin{equation*}
        \P\big(I(t)\not\subset t[-1/a,1/a]^d\big) \leq \exp\left(-t(1/a-2d\e n)\right) \frac{\e^{1+2den}}{\e-1}.
    \end{equation*}
    }
\end{cor}
\begin{proof}
    {\BBB We show $\P\big(I(at)\not\subset [-t,t]^d\big) \leq \e^{-t} (\e-1)^{-1}\e^{1+2den} \exp(2d\e at n)$ as the statement follows from the change of variables $\tilde{t}\mapsto at$.
    If $I(at)\not\subset [-t,t]^d$, then there exists some path $\gamma\colon o\leadsto x$ with $T^{(n)}(\gamma)\leq at$ and $|x|_\infty \geq t$. In particular, $\gamma$ starts in $o$ and has length at least $|\gamma|\geq|x|_\infty\geq t$. The probability that such a path exists can be estimated exactly as in (\ref{eq:prob-path-too-fast-1}) and (\ref{eq:prob-path-too-fast-2}), yielding the claim.}
\end{proof}

\begin{lemma}[Properties of the inverse speed]\label{lem:properties-inverse-speed} Let $x,y\in\R^d$ and $c\in\R$. Then,
    \begin{enumerate}[(i)]
        \item $\varphi$ is invariant under the symmetries of $\Z^d$ that fix the origin,
        \item $\varphi(x+y)\leq \varphi(x)+\varphi(y)$ and particularly $|\varphi(x)-\varphi(y)|\leq \varphi(x-y)$,
        \item $\varphi$ is continuous,
        \item $\varphi(cx)=|c| \varphi(x)$,
        \item $(2dn\e)^{-1}|x|_\infty \leq \varphi(x) \leq |x|_1 \leq d|x|_\infty$,
        \item $\varphi$ is Lipschitz continuous,
        \item \(\varphi\) does not depend on the {\BBB chosen atom-free distribution governing  the contact-times.}
    \end{enumerate}
\end{lemma}

\begin{proof}
Part~$(i)$ follows from the fact that the function $\varphi$ inherits the symmetry of the underlying model on $\Z^d$.

In order to prove subadditivity in Part~$(ii)$, recall 
\[
	\varphi(x)=\lim_{t\uparrow\infty}t^{-1}D(0,[tx])=\lim_{t\uparrow\infty}t^{-1}\widetilde{D}_0(0,[tx]).
\]
We further obtain 
\[
 	t^{-1}\widetilde{D}_0(0,[t(x+y)]) \leq t^{-1}\widetilde{D}_0(0,[tx])+t^{-1}\widetilde{D}_{\widetilde{D}_0(0,[tx])}([tx],[t(x+y)]).
 \]
Clearly, the left-hand side converges almost surely to $\varphi(x+y)$ and the first term on the right-hand side converges almost surely to $\varphi(x)$. Therefore, it remains to show 
\[
	\liminf_{t\uparrow\infty} t^{-1}\widetilde{D}_{\widetilde{D}_0(0,[tx])}([tx],[t(x+y)]) \leq \varphi(y).
\]
For the second term on the right-hand side, we use shift-invariance on the lattice and shift-invariance over days to see that the term is equal in law to
\[
	t^{-1}\widetilde{D}_0(0,[t(x+y)]-[tx]).
\]
We aim to compare this term to \(\widetilde{D}_0(0,[ty])/t\), which converges to $\varphi(y)$ almost surely. To this end, observe that
\[ 
	 t^{-1}\big|\widetilde{D}_0(0,[t(x+y)]-[tx])-\widetilde{D}_0(0,[ty])\big| \xrightarrow[]{} 0,
\]
in distribution, as \(t\to\infty\), since $|([t(x+y)]-[tx])-[ty]|\leq 3$. In particular, this implies
\[
	t^{-1}\widetilde{D}_{\widetilde{D}_0(0,[tx])}([tx],[t(x+y)]) \xrightarrow[]{} \varphi(y),
\]
in distribution. Since $\varphi(y)$ is a deterministic limit, the established convergence also holds in probability and we can therefore find a subsequence along which the convergence holds almost surely. This finishes the proof of subadditivity. The second statement in Part~$(ii)$ simply follows from symmetry and the reverse triangle inequality. 

The continuity of $\varphi$ is a result of Lemma \ref{lem:existence-inverse-speed-general-direction}, proving Part~$(iii)$.

In order to prove Part~$(iv)$, observe that we have, for all $s\in\N$,  
\[
	\varphi(sx)=\lim\limits_{t\uparrow\infty}t^{-1}D(0,tsx)
    	= \lim\limits_{t\uparrow \infty} s t^{-1}D(0,tx)  
    	= s\varphi(x)
\] 
and by distributional symmetry clearly $\varphi(-x)=\varphi(x)$.
Further, for each rational $c=p/m$
\[
	m\varphi(xp/m)= \varphi(px)=p\varphi(x),
\]
implying $\varphi(cx)=c\varphi(x)$. The claim for arbitrary $c\in\R$ then follows from the continuity in Part~$(iii)$. 

The upper bound in Part~$(v)$ is trivial since at least one edge per day is traversed. Furthermore, the shortest path from the origin $o$ to any $[kx]$ has length at most $|[kx]|_1\leq d(k|x|_\infty+1)$ and the lower bound is thus a consequence of Lemma \ref{lem:finite-speed-general-n} and {\BBB Part~$(iv)$ as $\varphi(x)=|x|_\infty \varphi(x/|x|_\infty)\geq (2dne)^{-1}|x|_\infty$.}

For the Lipschitz-continuity, we first observe that {\BBB by $(ii)$ and $(i)$
\[\sup_{|v|_1=1}\varphi(v) \leq \sup_{|v|_1=1}\sum_{i=1}^d |e_i\cdot v|\varphi(e_i) = \sup_{|v|_1=1}\sum_{i=1}^d |e_i\cdot v|\varphi(e_1) = \varphi(e_1)\,.\] 

Hence, for any $x\neq y$ by~$(ii)$ and~$(iv)$,}
\[
	\big|\varphi(x)-\varphi(y)\big|
		\leq\varphi(x-y)
		\leq |x-y|_1 \varphi\left(\frac{x-y}{|x-y|_1}\right) 
		\leq |x-y|_1 \cdot \varphi(e_1),
\]

and thus proving Part~$(vi)$.

Finally, the independence of \(\varphi\) {\BBB on the underlying chosen atom-free distribution of the $U^{(k,i)}$ follows from the observation that \(\varphi\) only depends on the order statistics of the contact times, which is the same for all atom-free distributions.} This concludes the proof. 
\end{proof}

We have now collected all the results required to prove \Cref{thm:shape}. Let us first define the limiting shape \(\Bcal=\Bcal_n\) via
\begin{equation}\label{eq:limiting_shape}
	\mathcal{B}:=\{x\in\R^d\colon \varphi(x)\leq 1\},
\end{equation}
and note that, due to Lemma \ref{lem:properties-inverse-speed}, $\mathcal{B}$ is indeed closed and convex, $\mathcal{B}\subset[-b,b]^d$ for $b=2dn\e$, and \(\mathcal{B}\) does not depend on the distribution of the contact times.

\begin{proof}[Proof of \Cref{thm:shape}]
We argue by contradiction. Assume that the theorem does not hold. That is, there exists some $\eps>0$ and an event \(\Omega_\varepsilon\) of positive probability, on which {\BBB we find a sequence $(t_k)_{k\in\N}\subset\R$ with $t_k\uparrow\infty$, together with $z_{t_k}\in\Z^d$} such that either
\begin{equation}\label{eq:case1-shape-theorem}
    z_{t_k}\in I(t) \text{ and } z_{t_k}\notin (1+\eps)t\mathcal{B}
\end{equation}
or
\begin{equation}\label{eq:case2-shape-theorem}
    z_{t_k}\notin I(t) \text{ and } z_{t_k}\in (1-\eps){t_k}\mathcal{B} \,.
\end{equation}
Let us focus on the case that~\eqref{eq:case2-shape-theorem} happens for infinitely many $k$. The case~\eqref{eq:case1-shape-theorem} is handled similarly as explained below.

As $\mathcal{B}$ is bounded, we find some $x\in\R^d$ and a  subsequence of $(t_k)_{k\in\N}$ {\BBB (for simplicity still denoted by $(t_k)_{k\in\N}$)} such that
\[
	z_{t_k}\in (1-\eps)t_k\mathcal{B}  \qquad\text{and}\qquad \lim_{k\uparrow\infty}z_{t_k} / t_k = x.
\]
Furthermore, since $\mathcal{B}$ is closed, we have $x\in(1-\eps)\mathcal{B}$. In particular, by Lemma~\ref{lem:existence-inverse-speed-general-direction} and Definition~\eqref{eq:limiting_shape}, we have 
\[
	\lim_{t\uparrow\infty} t^{-1}D_0(0,[tx])=\varphi(x) \leq 1-\eps.
\]
However, $z_{t_k}\notin \tilde{I}(t_k)$ is equivalent to $D_0(0,z_{t_k})>t_k$, i.e.,
\[
	D_0(0,z_{t_k})/t_k > 1,
\]
 which yields, by use of the reverse triangle inequality~\eqref{eq:reverse-triangle-real-points} for $z_{t_k}$ and $t_k x$,
\begin{equation*}
    1 <  D_0(0,z_{t_k})/t_k 
      \leq |z_{t_k}/t_k - x|_1 + 2d/t_k +  D_0(0,[t_kx])/t_k 
      \xrightarrow[]{k\to\infty} \varphi(x)  \leq 1-\eps,
\end{equation*}
providing a contradiction. Let us briefly discuss the remaining case, and assume that~\eqref{eq:case1-shape-theorem} happens infinitely often. For some fixed $a<2dn\e$, Corollary~\ref{cor:exp-decay-general-n} together with the Borel--Cantelli Lemma implies  $\tilde{I}(t)/t\subset [-1/a,1/a]^d$ almost surely for large $t$. We infer the existence of a convergent subsequence $z_{t_k} / t_k\xrightarrow[]{} x$, as \(k\to\infty\). Analogous calculations yield the contradiction
\[
	1 \geq D_0(0,z_{t_k})/t_k \geq D_0(0,[t_kx])/t_k  
		-\big( |z_{t_k}/t_k - x|_1 + 2d/t_k\big) 
      \xrightarrow[]{k\to\infty} \varphi(x) \geq 1+\eps.
\]
This concludes the proof.
\end{proof}

\subsection{Proof of couplings}\label{sec:proof_couple}
We provide various couplings of the model based on~\eqref{eq:model_n_Day} in this section. This will ultimately result in showing Proposition \ref{prop:speed-up-through-rigidity} and \Cref{thm:approx}, i.e., verifying that the rescaled shape \(\Bcal_n/n\) converges w.r.t.\ the Hausdorff metric towards the limiting shape \(\Bcal_R\) of the Richardson model.

We start with the proof of Proposition~\ref{prop:speed-up-through-rigidity}, which enables a simpler proof of \Cref{thm:approx} later. Loosely speaking, the rigid model~\eqref{eq:modelRigid} is faster since not being able to traverse an edge in a given day will give a bias to traverse said edge earlier the next day. {\BBB This will be made precise in the proof using the following strategy: The exploration processes of the infected sets associated to~\eqref{eq:model_n_Day} and~\eqref{eq:modelRigid}, respectively, are modeled jointly as a two-step random experiment, described by two independent uniform random variables assigned to each edge. When the respective infection first reaches an edge, the remaining number of future contacts on the day of discovery is jointly governed by the first uniform random variable. If there is at least one such contact, the exploration progresses by an edge traversal at the time of the first of the remaining contacts, which is encoded by the second uniform random variable. However, if the first step reveals that there is no remaining contact during the same day, then, in the~\eqref{eq:model_n_Day} model, the traversal happens uniformly at random on the next day, whereas, in the~\eqref{eq:modelRigid} model, the traversal time is still distributed uniformly but on a smaller interval. The length of this smaller interval is determined by the time of infection on the previous day. Notably, these distributions are both described by the second uniform random variable. In that way, the distributions of both explorations are governed by an i.i.d.\ family of two-dimensional uniforms and thus coupled. }

\begin{proof}[Proof of Proposition \ref{prop:speed-up-through-rigidity}]
	Let \((U^{(1)}_e)_{e\in E}\) and \((U_e^{(2)})_{e\in E}\) be two independent sequences of independent uniformly on \((0,1)\) distributed random variables. We start by describing the exploration processes corresponding to the models~\eqref{eq:model_n_Day} and~\eqref{eq:modelRigid} in terms of those uniform random variables, thus coupling the two explorations. To this end, denote by
	\[
		F_{n,p}(\cdot)=\sum_{i=0}^{\lfloor\cdot\rfloor} \binom{n}{i}p^i(1-p)^{n-i}
	\]
	the distribution function of a binomially distributed random variable. For each edge \(e\), we define the function \(K_e: [0,1]\to\{0,\dots,n\}\) as  
	\begin{equation*}\label{eq:def-K-in-coupling}
        K_e(t)=\min \big\{ k\in\{0,\dots,n\} \colon U^{(1)}_e \leq F_{n,(1-t)}(k) \big\},
    \end{equation*}
	and extend \(K_e\) to the whole real line by way of \(K_e(t):=K_e(t-\lfloor t\rfloor)\). The quantity \(K_e(t)\) thus models the number of contacts left on the edge \(e\) during \([t,\lceil t\rceil)\), where \(\lceil t \rceil = \lfloor t\rfloor +1\). Put differently, if the infection arrives at an edge for the first time on the \(\lfloor t\rfloor\)-th day at time \(t-\lfloor t\rfloor\), then there are \(K_e(t)\) many contacts still happening until the end of that day. By definition, it is easy to see that for \(0<s<t<1\), we have \(K_e(s)\geq K_e(t)\) since  \(1-t<1-s\). 
	
	For \(\tau>0\) and \(K\in\N\), we additionally define  the random variable
	\begin{equation*}\label{eq:def-T-in-coupling}
        T_e(\tau,K)=\tau\Big(1- \sqrt[K]{U^{(2)}_e}\Big),
    \end{equation*}
    noting that it equals, in distribution, the minimum of \(K\) independent random variables distributed uniformly on \([0,\tau]\).
    
    We proceed by describing the exploration process \(I\) associated to~\eqref{eq:model_n_Day} in terms of \(K_e\) and \(T_e\). Recall that, whenever the exploration discovers a new vertex that has not yet been seen at some time $s$, its edges incident to yet unexplored vertices are revealed, {\BBB meaning that we reveal all contact times from this point onwards}. 
    Consider such a fixed edge $e$. If there is at least one more contact happening during \([s,\lceil s\rceil]\), then the other end-vertex becomes revealed at the minimum of the remaining contact times. If there is no more contact happening, the other end-vertex becomes revealed at the minimum of the \(n\) newly sampled contact times within \([\lceil s\rceil, \lceil s\rceil +1]\). However, by the above, the probability of seeing at least one more contact during the day of discovery coincides with the probability of \(\{K_e(s)\geq 1\}\). Moreover, conditionally on \(\{K_e(s)\geq 1\}\), the next contact has the same distribution as
    \begin{equation}\label{eq:coupNextContact}
    	s+T_e\big(\lceil s\rceil -s, K_e(s)\big) = \lceil s\rceil -\big(\lceil s\rceil -s\big)\sqrt[K_e(s)]{U_e^{(2)}}.
    \end{equation}  
    Similarly, if there is no more contact happening on the same day, meaning \(K_e(s)=0\), the first contact on the \(\lceil s\rceil\)-th day has the same distribution as
    \begin{equation}\label{eq:coupNextDayInd}
    	\lceil s\rceil + T_e(1,n) = \lceil s\rceil +1 - \sqrt[n]{U_e^{(2)}}.
    \end{equation}
    
    Analogously, we describe \(I'\), the exploration process associated to~\eqref{eq:modelRigid}. First of all, if an edge \(e\) is considered for the first time at time \(s\) and there is still a contact {\BBB happening within \([s,\lceil s\rceil]\), this contact can again be described as in~\eqref{eq:coupNextContact}}. However, conditionally on \(K_e(s)=0\), the information that no contact falls into \([s,\lceil s\rceil]\) yields that the first contact time on the next day equals in distribution  
    \begin{equation}\label{eq:coupNextDayRig}
    	\lceil s\rceil + T_e\big(s-\lfloor s\rfloor ,n\big) = s+1- \big(s-\lfloor s\rfloor\big)\sqrt[n]{U_e^{(2)}}.
    \end{equation} 
    
    Having established the coupling {\BBB by way of using the same sequence of uniforms to describe the respective contact times}, it remains to prove its monotonicity in the sense that \(I(t)\subset I'(t)\) for all \(t\geq 0\). We work on the probability space where \(I\) and \(I'\) are jointly defined by the previous coupling. Let \((u^{(1)}_e)_e\) and \((u^{(2)}_e)\) be realizations of the used sequences of independent random variables. Let \((s_v)_{v\in \Z^d}\) and \((s'_v)_{v\in \Z^d}\) be the collection of first times each vertex is explored by the processes \(I\) and \(I'\) respectively. Let \(0=\sigma'_0<\sigma'_1,\dots\) be the elements of \((s'_v)_v\) sorted from smallest to largest value. We now show inductively that for all times \(\sigma'_n\), we have \(I(\sigma'_n)\subset I'(\sigma'_n)\), and all edges that have been revealed but not yet been traversed in \(I'\) are traversed in \(I'\) no later than in \(I\). This particularly implies \(I(t)\subset I'(t)\) for all \(t\geq 0\). We start with the \emph{base case} \(\sigma'_0=0\). At time \(0\), we have \(I(0)=I'(0)=\{o\}\), i.e., only the origin has been discovered yet. Moreover, all incident edges and their next contact time become revealed. By the coupling, the next contact time on each of these edges coincide.
    
    Now assume for the \emph{induction step} that our assumption holds for \(\sigma'_{n-1}\). We note that \(I(\sigma'_{n-1})\subset I'(\sigma'_{n-1})\) implies that all edges revealed but not traversed in \(I\) are either already traversed or at least revealed in \(I'\) as well. The second assumption then guarantees that the edge traversed at time \(\sigma_n'\) by \(I'\) has not been traversed earlier by \(I\) and thus \(I(\sigma_n')\subset I'(\sigma_n')\) still. Let \(v\) be the vertex explored at time \(\sigma_n'\), i.e., \(s'_v=\sigma_n'\) and therefore \(s_v\geq s'_v\). Consider the edges incident to \(v\) that have not yet been revealed and we denote by \(t_{e}\), \(t_e'\) the time of the first contact on \(e\) after being revealed by the respective explorations. To finish the proof of the induction step, we have to show that \(t_e\geq t_e'\) for all newly revealed edges. To this end, we make use of~\eqref{eq:coupNextContact},~\eqref{eq:coupNextDayInd}, and~\eqref{eq:coupNextDayRig} and distinguish the following cases:
    \begin{enumerate}[(i)]
    	\item 
    		If \(\lceil s_v\rceil = \lceil s_v'\rceil\) and \(K_e(s_v)\geq 1\), then
    		\[
    			t_e' = \lceil s_v'\rceil-\big(\lceil s_v'\rceil-s_v'\big)\sqrt[K_e(s_v')]{u_e^{(2)}}\leq \lceil s_v\rceil-\big(\lceil s_v\rceil-s_v\big)\sqrt[K_e(s_v)]{u_e^{(2)}}=t_e,
    		\]
    		using \(K_e(s_v')\geq K_e(s_v)\) here, as outlined above. 
    	\item 
    		If \(\lceil s_v\rceil = \lceil s_v'\rceil\) and \(K_e(s_v')>K_e(s_v)=0\), then \(t_e'\leq\lceil s_v\rceil\leq t_e\).
    	\item
    		If \(\lceil s_v\rceil = \lceil s_v'\rceil\) and \(K_e(s_v')=K_e(s_v)=0\), then
    		\[
    			t_e'=\lceil s_v\rceil + T_e(s_v'-\lceil s'_v\rceil,n) \leq \lceil s_v\rceil + T_e(1,n) = t_e.
    		\]
    	\item
    		In all remaining cases, we have \(\lceil s'_v\rceil<s_v\). If additionally either \(s_v\geq \lceil s'_v\rceil+1\) or \(K_e(s_v)=0\), then \(t_e'\leq \lceil s'_v\rceil+1\leq t_e\). Similarly, if \(K_e(s_v')>0\), then \(t_e'\leq \lceil s_v'\rceil\leq s_v \leq t_e\). Finally, if \(K_e(s_v')=0\), \(K_e(s_v)>0\), and \(\lceil s'_v\rceil=\lfloor s_v\rfloor\), then \(\lceil s_v'\rceil+1= \lceil s_v \rceil\) and thus
            \begin{equation*}
                \begin{aligned}
                    t_e'&= s_v'+1 - \big(s_v'-\lfloor s_v'\rfloor\big)\sqrt[n]{u_e^{(2)}}
                \leq 
                \lceil s_v'\rceil + 1 - \big(s_v'-\lfloor s_v'\rfloor\big)\sqrt[K_e(s_v)]{u_e^{(2)}} - (\lceil s_v' \rceil -s_v')\sqrt[K_e(s_v)]{u_e^{(2)}} \\
                &= \lceil s_v\rceil - \big(\lceil s_v'\rceil - \lfloor s_v'\rfloor\big)\sqrt[K_e(s_v)]{u_e^{(2)}}
                \leq
                \lceil s_v\rceil -\big(\lceil s_v\rceil -s_v\big)\sqrt[K_e(s_v)]{u_e^{(2)}} =t_e.
                \end{aligned}
            \end{equation*}
    \end{enumerate}   
    Summarizing, we have \(t_e'\leq t_e\) for all newly revealed edges, which finishes the induction and thus the proof of Proposition~\ref{prop:speed-up-through-rigidity}. 
\end{proof}

Now, let us move towards the proof of \Cref{thm:approx}. We begin with the following observations. First, we may restrict ourselves without loss of generality to the special case of the model based on~\eqref{eq:model_n_Day}, in which the contact times are distributed uniformly on \([0,1)\) as this does not influence the limiting shape \(\Bcal_n\). Put differently, the random variables \(\{U_e^{(k,i)}\colon k\in\N_0,i=1,\dots,n\}\) on each edge form i.i.d.\ sequences of \(\operatorname{Uniform}(0,1)\) random variables. Secondly, we want to slow down time by a factor of \(n\). Instead of assigning each edge \(n\) contact times per day (an interval of unit length), we assign the \(n\) contact times to a period of \(n\) days represented by an interval of the form \([kn,(k+1)n)\). Formally, assign independently to each edge \(e\) the independent random variables \(\{\overline{U}_e^{(k,i)}\colon k\in\N_0,\, i=1,\dots,n\}\) distributed uniformly on \([kn,(k+1)n)\).
Then, the corresponding FCP is defined via the point processes 
\begin{equation}\label{eq:modelDay_indep}  
	\overline{\X}^{(n)}_e := \bigcup_{k=0}^\infty \bigcup_{i=1}^n \big\{\overline{U}_e^{(k,i)}\big\}.
\end{equation}
Let us denote the corresponding process of reachable vertices and the limiting shape by \(\overline{I}(t)\), respectively \(\overline{\Bcal}_n\), which are equal to \(I(t/n)\), respectively \(\Bcal_n/n\).  

In order to prove \Cref{thm:approx}, we rely on the following continuity result~\cite[Theorem~3]{cox1981continuityFPP}.
\begin{thm}[{Continuity of the time constant in FPP~\cite[Theorem~3]{cox1981continuityFPP}}]\label{thm:continuity-FPP}
Let $\mu_n$ and $\mu$ be distributions supported on $[0,\infty)$ such that $\mu_n\xrightarrow{w} \mu$.
Then, for all $x\in\Z^d$
\begin{equation*}
   \lim_{n\uparrow\infty} \psi^{(\mu_n)}(x)=\psi^{(\mu)}(x),
\end{equation*}
where \(\psi^{(\mu)}(x)\) denotes the time constant in direction \(x\in\R^d\), meaning 
\[
	\psi^{(\mu)}(x) = \lim_{t\uparrow\infty} t^{-1} D^{(\mu)}(o,[tx]). 
\] 
\end{thm}
\begin{remark}
	The analogous continuity of the inverse speed as established for FPP in \Cref{thm:continuity-FPP} does not hold in general for FCP: Consider a version of the FCP model in which the contact times  are distributed uniformly on $[0,1/k]$. By \Cref{thm:shape} we observe the same inverse speed and the same limiting shape for each fixed \(k\). However, as $k\to\infty$, the limiting distribution has an atom at zero, which results in instant transmission across all of $\Z^d$.
\end{remark}

\begin{proof}[Proof of \Cref{thm:approx}] Let us first construct the coupling between the time-rescaled model of the model based on~\eqref{eq:model_n_Day} (i.e., having \(n\) contacts in \(n\) days) and the Richardson model. For an edge $e$ and time $t\geq 0$, let $\overline{X}^{(n)}_{e,t}$ be the time duration until the next contact from $t$ on the edge \(e\). That is, 
    \[
     	\overline{X}^{(n)}_{e,t}=\min\big\{s\geq 0\colon  t+s\in \overline{\X}^{(n)}_e \big\}.
    \]
    Given some \(t\in[0,n)\), we have for any \(s\in[0,n)\) with \(s+t\leq n\), 
    \begin{equation}\label{eq:FCP-waiting-time-distributionI}
        \P(\overline{X}^{(n)}_{e,t} > s) = \P\big(\overline{\X}^{(n)}_e\cap [t,t+s] = \emptyset\big) = \left(1-s/n\right)^n,
    \end{equation}
    while we have for $s\in[0,2n]$ with $2n \geq t+s\geq n$, by the independence of the contact times in the intervals \([0,n)\) and \([n,2n]\),
    \begin{equation}\label{eq:FCP-waiting-time-distributionII}
    	\begin{aligned}
    		\P(\overline{X}^{(n)}_{e,t} > s) 
    			&
    				= \P(\overline{\X}^{(n)}_e\cap [t,t+s] = \emptyset)
               \\ &
               		= \P\big(\overline{\X}^{(n)}_e\cap [t,n) = \emptyset\big) \, \P\big(\overline{\X}^{(n)}_e\cap [n,t+s) = \emptyset\big)
               \\ &
               		= \Big(\frac{t}{n}\Big)^n  \Big(\frac{2n-(t+s)}{n}\Big)^n.
    	\end{aligned}
    \end{equation}
    Note that this is true for any \(t\in[kn,(k+1)n)\) for some \(k\) by shift invariance. Furthermore, note that~\eqref{eq:FCP-waiting-time-distributionI} is included in~\eqref{eq:FCP-waiting-time-distributionII} by choosing $t=n$.
    We observe from Equation~\eqref{eq:FCP-waiting-time-distributionII} that
    \begin{equation*}
        \begin{aligned}
            \P\big(\overline{X}^{(n)}_{e,t} > s\big) &= \Big(\frac{t}{n}\Big)^n  \Big(\frac{2n-(t+s)}{n}\Big)^n 
             = \Big(1+\frac{t-n}{n}\Big)^n  \Big(1+\frac{n-(t+s)}{n}\Big)^n
            \leq \e^{t-n+n-(t+s)} = \e^{-s},
        \end{aligned}
    \end{equation*}
    which shows the claim \(\overline{I}(t)\supset J_R(t)\) via a standard coupling in the exploration process.
    
    Next, we show the convergence of \(\overline{\Bcal}_n/n\). Using, Proposition \ref{prop:speed-up-through-rigidity}, it suffices to show \(\Bcal_n'/n\to\Bcal_R\).
    As the latter is a FCP model given by an underlying stationary random closed set, Proposition \ref{prop:shape-stationary} allows us to equivalently consider an FPP model with waiting time distribution \(\mu^{(n)}\)given by~\eqref{eq:FCP-FPP-link}. We see that for arbitrary \(s\geq0\) (and \(n\geq s\)), we have
    \begin{equation*}
        \mu^{(n)}([s,\infty))=(1-s/n)^n\xrightarrow{n\to\infty}\e^{-s}.
    \end{equation*}
    Since the (deterministic) limiting shapes are solely defined via the (deterministic) time constants, \Cref{thm:continuity-FPP} finishes the proof.
\end{proof}

\paragraph{Acknowledgement.} The authors thank Julia H\"{o}rrmann for her contributions to the project and thank her as well as Christian Hirsch, Jonas K\"oppl, Jean-Baptiste Gou\'er\'e, and Daniel Ahlberg for inspiring discussions. This research was supported by the Leibniz Association within the Leibniz Junior Research Group on {\em Probabilistic Methods for Dynamic Communication Networks} as part of the Leibniz Competition (grant no.\ J105/2020) and the Berlin Cluster of Excellence {\em MATH+} through the project {\em EF45-3} on {\em Data Transmission in Dynamical Random Networks}.

\section*{References}
\renewcommand*{\bibfont}{\footnotesize}
\printbibliography[heading = none]
\end{document}